\documentclass[11pt,a4paper,reqno]{amsart}
\usepackage[latin1]{inputenc}
\usepackage{amsmath}
\usepackage{amsfonts}
\usepackage{mathrsfs}
\usepackage{amssymb}
\usepackage{amsthm}
\usepackage{anysize}
\usepackage{enumitem}
\usepackage{amscd}
\usepackage{hyperref}
\usepackage[square, comma, numbers, sort&compress]{natbib}
\usepackage{indentfirst}
\marginsize{3.0cm}{3.0cm}{3.0cm}{3.0cm}
\setenumerate{label={\normalfont(\arabic*)}}
\newcommand{\form}[1]{{\langle #1 \rangle }}

\newcommand{\pfister}[1]{{\langle \! \langle #1 \rangle \! \rangle}}

\newcommand{\mydim}[1]{{\mathrm{dim}\!\; #1}}

\newcommand{\anispart}[1]{#1_{\mathrm{an}}}
\newcommand{\windex}[1]{{\mathfrak{i}_W(#1)}}
\newcommand{\witti}[2]{{\mathfrak{i}_{#1}(#2)}}

\newcommand{\simform}[1]{{#1_{\mathrm{sim}}}}
\newcommand{\normform}[1]{{#1_{\mathrm{nor}}}}

\newtheorem{theorem}{Theorem}[section]

\newtheorem{lemma}[theorem]{Lemma}
\newtheorem{claim}[theorem]{Claim}

\newtheorem{proposition}[theorem]{Proposition}
\newtheorem{corollary}[theorem]{Corollary}
\newtheorem{conjecture}[theorem]{Conjecture}

\theoremstyle{definition}

\theoremstyle{remark}
\newtheorem{remark}[theorem]{Remark}
\newtheorem{remarks}[theorem]{Remarks}

\numberwithin{equation}{section}
\begin{document}

\title[On the holes in $I^n$ in characteristic 2]{On the holes in $I^n$ for symmetric bilinear forms in characteristic 2}
\author{Stephen Scully}
\address{Department of Mathematics and Statistics, University of Victoria}
\email{scully@uvic.ca}


\maketitle

\begin{abstract} Let $F$ be a field. Following the resolution of Milnor's conjecture relating the graded Witt ring of $F$ to its mod-2 Milnor $K$-theory, a major problem in the theory of symmetric bilinear forms is to understand, for any positive integer $n$, the low-dimensional part of $I^n(F)$, the $n$th power of the fundamental ideal in the Witt ring of $F$. In a 2004 paper, Karpenko used methods from the theory of algebraic cycles to show that if $\mathfrak{b}$ is a non-zero anisotropic symmetric bilinear form of dimension $< 2^{n+1}$ representing an element of $I^n(F)$, then $\mathfrak{b}$ has dimension $2^{n+1} - 2^i$ for some $1 \leq i \leq n$. When $i = n$, a classical result of Arason and Pfister says that $\mathfrak{b}$ is similar to an $n$-fold Pfister form. At the next level, it has been conjectured that if $n \geq 2$ and $i= n-1$, then $\mathfrak{b}$ is isometric to the tensor product of an $(n-2)$-fold Pfister form and a $6$-dimensional form of trivial discriminant. This has only been shown to be true, however, when $n = 2$, or when $n = 3$ and $\mathrm{char}(F) \neq 2$ (another result of Pfister). In the present article, we prove the conjecture for all values of $n$ in the case where $\mathrm{char}(F) =2$. In addition, we give a short and elementary proof of Karpenko's theorem in the characteristic-2 case, rendering it free from the use of subtle algebraic-geometric tools. Finally, we consider the question of whether additional dimension gaps can appear among the anisotropic forms of dimension $\geq 2^{n+1}$ representing an element of $I^n(F)$. When $\mathrm{char}(F) \neq 2$, a result of Vishik asserts that there are no such gaps, but the situation seems to be less clear when $\mathrm{char}(F) = 2$.
 \end{abstract}
 
\section{Introduction} \label{SECintroduction} 

Let $F$ be a field. Recall that the Witt ring of symmetric bilinear forms over $F$, $W(F)$, is defined by taking the Grothendieck completion of the semiring of isometry classes of non-degenerate symmetric bilinear forms over $F$ (under the operations of orthogonal sum and tensor product) and then quotienting out the ideal generated by the hyperbolic plane $\mathbb{H}$. Each element of $W(F)$ is represented by the class of an anisotropic form which is unique up to isometry, and the classes of the forms of even dimension constitute a maximal ideal known as the fundamental ideal and denoted $I(F)$. By Witt decomposition, the problem of classifying all symmetric bilinear forms over $F$ amounts to the determination of $W(F)$, and the $I(F)$-adic filtration of the latter has played a central role in the subject since it was first investigated in the 1960s by Pfister. Indeed, a key driver in the overall development of the subject in the decades that followed was the now-famous conjecture of Milnor proposing a canonical identification of the graded ring $\mathrm{gr}_{I^\bullet(F)}W(F)$ with the mod-2 Milnor $K$-theory ring $K_\bullet^M(F)/2$. In the 1982 article \cite{Kato}, Kato gave a relatively direct proof of Milnor's conjecture in the case where $\mathrm{char}(F) = 2$ by making a fundamental connection between the two graded rings of interest and the differential graded algebra of absolute K\"{a}hler differentials of $F$ (in that setting). The case where $\mathrm{char}(F) \neq 2$, on the other hand, seemingly demanded the use of much more sophisticated algebraic-geometric tools, and was not resolved until the mid 1990s in celebrated work of Voevodsky (\cite{Voevodsky}) and Orlov-Vishik-Voevodsky (\cite{OVV}). In the characteristic-not-2 setting, Milnor's conjecture establishes a deep link between symmetric bilinear forms and motivic cohomology with $\mathbb{Z}/2\mathbb{Z}$ coefficients, and this was further enhanced by work of Morel (\cite{Morel1, Morel2}), who showed that the fibred products $I^n(F) \times_{K_n^M(F)/2}K_n^M(F)$ are canonically identified with certain stable homotopy groups of the motivic sphere spectrum of $F$. As such, the study of the Witt ring and the filtration induced by the fundamental ideal may be perceived as something that pertains to the very foundation of motivic algebraic geometry.

Following the resolution of Milnor's conjecture, one of the major outstanding problems in the theory of symmetric bilinear forms is to understand the ``low-dimensional'' part of the ideal $I^n(F)$ for each positive integer $n \geq 2$. More precisely, it is expected here that there is a reasonable classification of the anisotropic forms of dimension strictly less than $2^{n+1}$ that represent an element of $I^n(F)$. Perhaps the most significant work underlying this expectation is a 2004 article of Karpenko (\cite{Karpenko2}) in which strong constraints on certain discrete invariants of the relevant forms were established. In particular, Karpenko proved the existence of non-trivial ``dimension gaps'' (or ``holes'') in $I^n(F)$:

\begin{theorem}[Dimension Gap Theorem, \cite{Karpenko2, Laghribi2}] \label{THMdimensiongaptheorem} Let $n$ be a positive integer, and let $\mathfrak{b}$ be an anisotropic symmetric bilinear form of dimension $<2^{n+1}$ over $F$. If $\mathfrak{b}$ represents an element of the ideal $I^n(F)$, then $\mydim{\mathfrak{b}} = 2^{n+1} - 2^i$ for some integer $1 \leq i \leq n+1$. 
\end{theorem}

In \cite{Karpenko2}, it was in fact assumed that $\mathrm{char}(F) \neq 2$. However, it was later noted by Laghribi in \cite{Laghribi2} that the characteristic-2 case of Theorem \ref{THMdimensiongaptheorem} may be deduced from the characteristic-0 case by a routine specialization argument, and so Karpenko's work ultimately covers all cases. The existence of the first gap from $0$ to $2^n$ is the so-called Hauptsatz of Arason-Pfister, already established in 1970 as part of the early investigations into Milnor's conjecture.\footnote{As pointed out by Detlev Hoffmann, the treatment of the characteristic-2 case of the Hauptsatz given in \cite{ArasonPfister} is flawed. A correct proof was given in \cite{Laghribi2}, however, and another can be found in \S \ref{SECgaps} below.} While the Hauptsatz was achieved by direct, elementary means, Karpenko's proof of Theorem \ref{THMdimensiongaptheorem} relies on methods from the theory of algebraic cycles that may be deemed non-elementary due to an essential and somewhat mysterious use of Steenrod operations for Chow groups with $\mathbb{Z}/2\mathbb{Z}$ coefficients. This built on work of Vishik (\cite{Vishik}), who used similar methods to establish the existence of the second gap from $2^n$ to $2^n + 2^{n-1}$. Vishik also showed that in characteristic not 2, Theorem \ref{THMdimensiongaptheorem} cannot be improved in general. More precisely, it is easy to see that the dimensions permitted by the statement can all be realized, and so the interesting part of Vishik's observation is that there are no dimension gaps outside the range under consideration: If $\mathrm{char}(F) \neq 2$, and $m$ is any even integer $\geq 2^{n+1}$, then there exists a field extension $K/F$ and an anisotropic form of dimension $m$ over $F$ that represents an element of $I^n(K)$. We remark that Vishik's proof of this result also relies on algebraic-geometric methods that require the ground field to have characteristic different from 2 (see \cite[Prop. 82.7]{EKM}), and, unlike the situation for Theorem \ref{THMdimensiongaptheorem} itself, optimality in characteristic 2 no longer follows in an obvious way.

As far as classifying the low-dimensional elements of the ideals $I^n(F)$ goes, only a partial vision seems to have emerged. At the first level, the standard proof of the Hauptsatz shows that the forms of dimension $2^n$ that represent an element of $I^n(F)$ are precisely the general $n$-fold bilinear Pfister forms, i.e., tensor products of $n$ binary forms. Prior to this, Pfister (\cite{Pfister}) also showed that in the case where $\mathrm{char}(F) \neq 2$, the $6$-dimensional forms representing an element of $I^2(F)$ are precisely the Albert forms (i.e., $6$-dimensional forms of trivial discriminant), and the $12$-dimensional forms representing an element of $I^3(F)$ are precisely those that decompose as the tensor product of a binary form and an Albert form. This suggests that, at the second level, the following may hold:

\begin{conjecture} \label{CONJAlbertconj} Let $n$ be an integer $\geq 2$, and let $\mathfrak{b}$ be an anisotropic symmetric bilinear form of dimension $2^n + 2^{n-1}$ over $F$. If $\mathfrak{b}$ represents an element of $I^n(F)$, then $\mathfrak{b} \simeq \mathfrak{c} \otimes \mathfrak{d}$ for some $(n-2)$-fold Pfister form $\mathfrak{c}$ and Albert form $\mathfrak{d}$ over $F$.
\end{conjecture}

Beyond Pfister's results, the main evidence for this conjecture comes in the case where $\mathrm{char}(F) \neq 2$ and $n=4$. Here, Hoffmann (\cite{Hoffmann2}) reduced the problem to showing that any 24-dimensional form representing an element of $I^4(F)$ splits over a quadratic extension of $F$. In the recent article \cite{HarveyKarpenko}, Harvey and Karpenko have shown that any such form splits over an extension of degree $2m$ for some odd integer $m$, but a complete solution is yet to be established. As for the low-dimensional elements of $I^n(F)$ that lie beyond the second level, a classification of the $14$-dimensional forms representing an element of $I^3(F)$ has been achieved by Rost in the case where $\mathrm{char}(F) \neq 2$ (\cite{Rost}), but we are unaware of any precise conjectures about the general classification in this range. 

Turning to the present work, while some aspects of the theory of symmetric bilinear forms become more complicated in characteristic $2$, others are less so. We have already indicated, for instance, the the characteristic-2 case of Milnor's conjecture turned out to be substantially simpler than the characteristic-not-2 case. One might therefore hope that the study of the low-dimensional elements in the ideals $I^n(F)$ is more approachable in the characteristic-2 setting. In this article, we make some contributions in this direction. First, we give in \S \ref{SECgaps} a direct and elementary proof of the characteristic-2 case of Theorem \ref{THMdimensiongaptheorem}. Recall that the (only) previously existing proof here relies on the characteristic-0 case, which was proved using subtle algebraic-geometric methods. By contrast, the proof given here is no more involved that the proof of the Hauptsatz, and could have been found long before Theorem \ref{THMdimensiongaptheorem} was first conjectured in \cite{Vishik}. Second, while we are generally unable to show that Vishik's result on the absence of additional dimension gaps in the characteristic-not-2 setting extends to characteristic 2, we can at least show that this is the case when $n \leq 5$, among some additional observations (see Proposition \ref{PROPmainoptimality} and Corollary \ref{CORoptimalityfornleq5}). Third, and most significantly, we show in \S \ref{SECmain} the following:

\begin{theorem} \label{THMmain} Conjecture \ref{CONJAlbertconj} holds in the case where $\mathrm{char}(F) = 2$. \end{theorem}

Unlike the situation in which $\mathrm{char}(F) \neq 2$, the $n=3$ case of this result already seems to have been absent from the literature. We would like to remark here that a proof of this particular case was privately communicated to the author by Detlev Hoffmann several years ago, but his approach relied on results of Faivre on $24$-dimensional non-singular quadratic forms (\cite{Faivre}), and cannot be adapted to the cases where $n \geq 4$ at the present time.\footnote{Effectively, one would first have to prove Conjecture \ref{CONJAlbertconj} in the case where $\mathrm{char}(F) \neq 2$ to make this approach work in general.} The key tool in our proof of Theorem \ref{THMmain} is \cite[Thm. 6.4]{Scully1}, which strongly constrains possible isotropy behaviour of quasilinear quadratic forms over function fields of quadrics (for far-reaching discrete constraints that result from this, see the recent article \cite{Scully2}). As part of the proof, we also solve the quasilinear version of an open problem on non-singular quadratic forms whose first higher Witt index is exactly a quarter of the dimension of the form (see Theorem \ref{THMi12^{n-2}} and Remark \ref{REMi12^{n-2}}). Similar ideas should be applicable to the study of the remaining low-dimensional elements in $I^n(F)$, as well as the symmetric bilinear forms of Knebusch height $2$ (of which the forms considered in Conjecture \ref{CONJAlbertconj} are examples) in the characteristic-2 setting. These problems will be considered in subsequent work. 

\section{Preliminaries} \label{SECprelims}

For the remainder of this article, $F$ will denote an arbitrary field of characteristic 2. In this section, we present some preliminary material on symmetric bilinear and quasilinear quadratic forms that will be needed later on. Short proofs are given where possible. For any notation, terminology or basic facts that are not addressed here, we refer to \cite{EKM}. 

\subsection{Symmetric Bilinear Forms} \label{SUBSECbilinearforms} By a \emph{symmetric bilinear form over $F$}, we will mean a pair $(V,\mathfrak{b})$ consisting of a finite-dimensional $F$-vector space $V$ and a non-degenerate symmetric $F$-bilinear map $\mathfrak{b} \colon V \times V \rightarrow F$. In practice, we supress $V$ from our notation and simply talk about the form $\mathfrak{b}$. The dimension of $\mathfrak{b}$ (which is just the dimension of $V$) will be denoted $\mydim{\mathfrak{b}}$. The quadratic form $V \rightarrow F$ that sends $v \in V$ to $b(v,v)$ will be denoted $\phi_{\mathfrak{b}}$. Since $\mathrm{char}(F) = 2$, $\phi_{\mathfrak{b}}$ is quasilinear in the sense of \S \ref{SUBSECquasilinear} below. By definition, $\mathfrak{b}$ is isotropic (or anisotropic) precisely when $\phi_{\mathfrak{b}}$ is. By a \emph{subform of $\mathfrak{b}$} we will mean any bilinear form isometric to an orthogonal summand of $\mathfrak{b}$. If $a_1,\hdots,a_n \in F^\times$, then we write $\form{a_1,\hdots,a_n}_b$ for the orthogonal sum $\perp_{i=1}^n \form{a_i}_b$, where $\form{a_i}_b \colon F \times F \rightarrow F$ is the $1$-dimensional form that sends $(x,y)$ to $a_ixy$. We also write $\pfister{a_1,\hdots,a_n}_b$ for the $n$-fold bilinear Pfister form $\form{1,a_1}_b \otimes \cdots \otimes \form{1,a_n}_b$. We recall that bilinear Pfister forms are \emph{round}, i.e., if $\mathfrak{b}$ is Pfister, then $\mathfrak{b} \simeq \lambda \mathfrak{b}$ for any $\lambda \in F^\times$ represented by the quadratic form $\phi_{\mathfrak{b}}$. A \emph{general} $n$-fold bilinear Pfister form is a form which is similar to (i.e., isometric to a non-zero scalar multiple of) an $n$-fold bilinear Pfister form. If $a \in F^\times$, then $\mathbb{M}_a$ denotes the metabolic plane $\form{a,a}_b$. Any isotropic symmetric bilinear form of dimension $2$ over $F$ is either isometric to some $\mathbb{M}_a$ or to the hyperbolic plane $\mathbb{H}$, which is up to isometry the unique symmetric bilinear form of dimension $2$ over $F$ whose associated quadratic form is zero. Witt decomposition asserts that any symmetric bilinear form $\mathfrak{b}$ over $F$ decomposes as an orthogonal sum of isotropic forms of dimension $2$ and an anisotropic form which, up to isometry, depends only on the isometry class of $\mathfrak{b}$. Using that $\mathbb{M}_a \perp \form{a}_b \simeq \mathbb{H} \perp \form{a}_b$ for any $a \in F^\times$, one can refine this as follows (see, e.g., \cite[(2.1)]{LaghribiMammone}):

\begin{proposition} \label{PROPWittdecomposition} Let $\mathfrak{b}$ be a symmetric bilinear form over $F$. Then there exist non-negative integers $\mathfrak{i}_{\mathbb{M}}(\mathfrak{b})$ and $\mathfrak{i}_{\mathbb{H}}(\mathfrak{b})$, an anisotropic symmetric bilinear form $\anispart{\mathfrak{b}}$ over $F$, and elements $a_1,\hdots,a_{\mathfrak{i}_{\mathbb{M}}(\mathfrak{b})} \in F^\times$ such that
\begin{itemize} \item $\mathfrak{b} \simeq \anispart{\mathfrak{b}} \perp \mathbb{M}_{a_1} \perp \cdots \perp \mathbb{M}_{a_{\mathfrak{i}_{\mathbb{M}}(\mathfrak{b})}} \perp \mathfrak{i}_{\mathbb{H}}(\mathfrak{b}) \cdot \mathbb{H}$, and
\item $\anispart{b} \perp \form{a_1,\hdots,a_{\mathfrak{i}_{\mathbb{M}}(\mathfrak{b})}}_b$ is anisotropic. \end{itemize}
Moreover, $\mathfrak{i}_{\mathbb{M}}(\mathfrak{b})$ and $\mathfrak{i}_{\mathbb{H}}(\mathfrak{b})$ are unique, and $\anispart{b}$ is unique up to isometry. 
\end{proposition}

In the situation of the proposition, the integer $\mathfrak{i}_{\mathbb{M}}(\mathfrak{b}) + \mathfrak{i}_{\mathbb{H}}(\mathfrak{b})$ is the Witt index of $\mathfrak{b}$ (i.e., the maximal dimension of a totally isotropic subspace for the quadratic form $\phi_{\mathfrak{b}}$), which we denote $\windex{\mathfrak{b}}$. We have the inequality $\windex{\mathfrak{b}} \leq \left[\frac{\mydim{\mathfrak{b}}}{2}\right]$, and we shall say that $\mathfrak{b}$ is \emph{split} when equality holds. For instance, any isotropic general Pfister form is split. The class of the form $\mathfrak{b}$ in the Witt ring $W(F)$ will be denoted $[\mathfrak{b}]$. There is a well-defined map $W(F) \rightarrow F^\times/F^{\times 2}$ taking $[\mathfrak{b}]$ to the determinant of $\mathfrak{b}$ (i.e., the square class of the determinant of any symmetric matrix representing $\mathfrak{b}$), which we denote $\mathrm{det}\;\mathfrak{b}$. This induces a group isomorphism $I(F)/I^2(F) \rightarrow F^\times/F^{\times 2}$ whose inverse takes the square class of an element $a \in F^\times$ to the class of the $1$-fold Pfister form $\pfister{a}_b$. In particular, the elements of $I^2(F)$ are precisely the classes of the anisotropic symmetric bilinear forms over $F$ which have trivial determinant. For any integer $n \geq 1$, it is clear that $I^n(F)$ is additively generated by the classes of the general $n$-fold Pfister forms over $F$. 

\subsection{Absolute K\"{a}hler Differentials and Kato's Theorem} \label{SUBSECKato} Let $(\Omega_F^\bullet,d)$ be the differential graded $F$-algebra of absolute K\"{a}hler differentials of $F$ (see \cite[Ch. 9]{Matsumura}). Recall that a subset $S$ of $F$ is said to be \emph{$2$-independent} if $[F^2(a_1,\hdots,a_m):F^2] = 2^m$ for every positive integer $m$ and every cardinality-$m$ set $\lbrace a_1,\hdots,a_m \rbrace \subset S$, and that a \emph{$2$-basis of $F$} is a maximal $2$-independent subset of $F$. By Zorn's Lemma, every $2$-independent subset of $F$ is contained in a $2$-basis of $F$. If $\mathcal{B}$ is a $2$-basis of $F$ with a strict total ordering $<$, then the set $\lbrace da_1 \wedge \cdots \wedge da_n\;|\; a_1,\hdots,a_n \in \mathcal{B}, a_1<a_2<\cdots <a_n \rbrace$ is a basis of $\Omega_F^n$ as an $F$-vector space for any non-negative integer $n$. As part of his work on the characteristic-2 case of Milnor's conjecture, Kato has shown the following fundamental result:

\begin{theorem}[\cite{Kato}] \label{THMKato} For any non-negative integer $n$, there exists an injective group homomorphism $e_n \colon I^n/I^{n+1} \rightarrow \Omega_F^n$ such that if $a_1,\hdots,a_n \in F^\times$, then $e_n$ maps the class of the $n$-fold bilinear Pfister form $\pfister{a_1,\hdots,a_n}_b$ to the $n$-fold wedge product $\frac{da_1}{a_1} \wedge \cdots \wedge \frac{da_n}{a_n}$. \end{theorem}

\begin{remark} For each non-negative integer $n$, there is a well-defined group homomorphism $\wp \colon \Omega_F^n \rightarrow \Omega_F^n/d(\Omega_F^{n-1})$ with the property that $\wp(a\frac{da_1}{a_1} \wedge \cdots \wedge \frac{da_n}{a_n}) = (a^2+a)\frac{da_1}{a_1} \wedge \cdots \wedge \frac{da_n}{a_n}$ for all $a,a_1,\hdots,a_n \in F^\times$. In \cite{Kato}, Kato has also shown that the image of the map $e_n$ in the statement of Theorem \ref{THMKato} is precisely the kernel of $\wp$. This will not be used below. \end{remark}

\subsection{Quasilinear Quadratic Forms} \label{SUBSECquasilinear} By a \emph{quasilinear quadratic form over $F$}, we will mean a pair $(V,\phi)$ consisting of a finite-dimensional $F$-vector space $V$ and a quadratic form $\phi \colon V \rightarrow F$ which is quasilinear, i.e., which has the additional property that $\phi(v + w) = \phi(v) + \phi(w)$ for all $v,w \in V$. As with symmetric bilinear forms, we supress $V$ from the notation and simply talk about the form $\phi$. The dimension of $\phi$ (which is just the dimension of $V$) will be denoted $\mydim{\phi}$. If $\mathfrak{b}$ is a symmetric bilinear form over $F$, then the quadratic form $\phi_{\mathfrak{b}}$ is quasilinear. Every quasilinear quadratic form over $F$ arises in this way, but (unlike the situation in characteristic different from $2$), the lift is very far from being unique (see \S \ref{SECbilinearlifts} below). If $a_1,\hdots,a_n \in F^\times$, then we write $\form{a_1,\hdots,a_n}$ for the quasilinear quadratic form $F^n \rightarrow F$ that sends $(x_1,\hdots,x_n)$ to $\sum_{i=1}^n a_ix_i^2$. In other words, $\form{a_1,\hdots,a_n} = \phi_{\form{a_1,\hdots,a_n}_b}$. Similarly, we write $\pfister{a_1,\hdots,a_n}$ for the quasilinear quadratic form associated to the $n$-fold bilinear Pfister form $\pfister{a_1,\hdots,a_n}_b$. Forms of this type are called \emph{$n$-fold quasi-Pfister forms}. Orthogonal sums and tensor products of quasilinear quadratic forms are defined in the obvious way, and are compatible with the analogous constructions for symmetric bilinear forms. If $\phi$ and $\psi$ are quasilinear quadratic forms over $F$, then we will say that $\psi$ is a \emph{subform} of $\phi$ (resp. that $\phi$ is \emph{divisible by $\psi$}) if $\phi \simeq \psi \perp \sigma$ (resp. $\phi \simeq \psi \otimes \sigma$) for some quasilinear quadratic form $\sigma$ over $F$. We use the notation $\psi \subset \phi$ to indicate that $\psi$ is a subform of $\phi$. If $\phi$ is a quasilinear quadratic form over $F$, then Witt decomposition asserts that $\phi \simeq \anispart{\phi} \perp \mathfrak{i}_0(\phi) \cdot \langle 0 \rangle$ for a unique non-negative integer $\witti{0}{\phi}$ and unique anisotropic quasilinear quadratic form $\anispart{\phi}$ over $F$ (up to isometry). We refer to the integer $\witti{0}{\phi}$ as the \emph{isotropy index of $\phi$}. When $\phi$ arises from a given symmetric bilinear form $\mathfrak{b}$, this relates to the Witt index of $\mathfrak{b}$ as follows:

\begin{lemma} \label{LEMWittindexrelation} If $\mathfrak{b}$ is a symmetric bilinear form over $F$, then $\witti{0}{\phi_{\mathfrak{b}}} = \windex{\mathfrak{b}} + \mathfrak{i}_{\mathbb{H}}(\mathfrak{b})$. 
\begin{proof} Immediate from Proposition \ref{PROPWittdecomposition}. 
\end{proof}
\end{lemma}

If $\phi$ is a quasilinear quadratic form over $F$, then we will write $D(\phi)$ for its value set in $F$. The additivity of $\phi$ implies that $D(\phi)$ is in fact an $F^2$-vector subspace of $F$. Moreover, up to isometry, $\anispart{\phi}$ is the unique anisotropic quasilinear quadratic form over $F$ whose value set coincides with that of $\phi$ (see \cite[Prop. 2.6]{Hoffmann3}). The form $\phi$ is therefore determined up to isometry by the pair $(D(\phi),\witti{0}{\phi})$, and if $\phi$ is anisotropic, then its subforms are precisely the anisotropic quasilinear quadratic forms over $F$ whose value sets are contained in $D(\phi)$. The anisotropic quasilinear quadratic forms over $F$ whose value sets are multiplicatively closed (and hence subfields of $F$) are precisely the anisotropic quasi-Pfister forms. This allows one to attach certain quasi-Pfister forms to any given quasilinear quadratic form $\phi$ over $F$. First, the subfield of $F$ generated by all products $xy$ with $x,y \in D(\phi)$ is called the \emph{norm field of $\phi$}, and is denoted $N(\phi)$. One readily observes that $N(\phi)$ is a finite purely inseparable extension of $F^2$, and we write $\mathrm{lndeg}(\phi)$ for the integer $\mathrm{log}_2([N(\phi):F^2])$, which we call the \emph{norm degree of $\phi$}.\footnote{In the literature, the term \emph{norm degree of $\phi$} usually refers to the $2$-power integer $[N(\phi):F^2]$, but we prefer to simply take the exponent of this $2$-power.} By the preceding remarks, there is, up to isometry, a unique anisotropic $\mathrm{lndeg}(\phi)$-fold quasi-Pfister form $\mathrm{lndeg}(\phi)$ over $F$ with $D(\normform{\phi}) = N(\phi)$. If $a \in D(\phi) \setminus \lbrace 0 \rbrace$, then $a\anispart{\phi} \subset \normform{\phi}$, and so $\mathrm{lndeg}(\phi) \geq s+1$, where $s$ is the unique integer satisfying $2^s < \mydim{\anispart{\phi}} \leq 2^{s+1}$. In the special case where $\phi$ is anisotropic and $\mathrm{lndeg}(\phi) = s+1$, we say that $\phi$ is a \emph{quasi-Pfister neighbour}. This is equivalent to saying that the (affine or projective) quadrics defined by $\phi$ and $\normform{\phi}$ are stably birationally equivalent as varieties over $F$ (see \cite[\S 2.G]{Scully2}). Returning to the general case, the set $G(\phi):= \lbrace a \in F\;|\; aD(\phi) \subseteq D(\phi) \rbrace$ is a subfield of $N(\phi)$ containing $F^2$ and is thus again the value set of a unique quasi-Pfister form $\simform{\phi} \subset \normform{\phi}$ over $F$ up to isometry. Concretely, $G(\phi)$ consists of 0 together with the similarity factors of $\phi$ in $F$, i.e., elements $a \in F^\times$ for which $a\phi \simeq \phi$. The set $D(\phi)$ is a $G(\phi)$-linear subspace of $F$, from which it follows that $\anispart{\phi}$ is divisible by $\simform{\phi}$. More generally, if $\pi$ is an anisotropic quasi-Pfister form over $F$, then $\anispart{\phi}$ is divisible by $\pi$ if and only if $\pi \subset \simform{\phi}$ (see, e.g., \cite[Lem. 2.11]{Scully2}). 

\subsection{Quasilinear Quadratic Forms and Field Extensions} \label{SUBSECfieldextensions} Let $\phi$ be a quasilinear quadratic form over $F$. For any field extension $L$ of $F$, the set $D(\phi_L)$ is nothing else but the $L^2$-linear span of $D(\phi)$ in $L$. Applying this to $\normform{\phi}$, we see that $N(\phi_L) = N(\anispart{(\phi_L)})$ is the compositum of $N(\phi)$ and $L^2$ inside $L$, and so $\normform{(\anispart{(\phi_L)})} \simeq \normform{(\phi_L)} \simeq \anispart{((\normform{\phi})_L)}$. We also have the following basic observations: 

\begin{lemma} \label{LEMbehaviouroverfieldextensions} Let $\phi$ be an anisotropic quasilinear quadratic form over $F$, and let $L$ be a field extension of $F$. Then:
\begin{enumerate}
\item If $\psi$ is a subform of $\phi$ that remains anisotropic over $L$, then $\psi_L \subset \anispart{(\phi_L)}$. 
\item If $\phi$ is divisible by a quasi-Pfister form $\pi$ over $F$, and $\pi$ remains anisotropic over $L$, then $\anispart{(\phi_L)}$ is divisible by $\pi_L$. \end{enumerate}
\begin{proof} (1) Clear from the above remarks.

(2) Let $a \in D(\pi) \setminus \lbrace 0 \rbrace$. Since $\phi$ is divisible by $\phi$, we have $ab \in D(\phi)$ for all $b \in D(\phi)$. Since $D(\phi_L)$ is the $L^2$-linear span of $D(\phi)$, it follows that $aD(\phi_L) \subseteq D(\phi_L)$, and so $a \in G(\anispart{(\phi_L)}) = D(\simform{(\anispart{(\phi_L)})})$. Since $D(\pi_L)$ is the $L^2$-linear span of $D(\pi)$, we then get that $D(\pi_L) \subseteq D(\simform{(\anispart{(\phi_L)})})$, and so $\pi_L \subset \simform{(\anispart{(\phi_L)})}$, proving the claim. \end{proof}
\end{lemma}

If $L$ is a separable field extension of $F$, then every anisotropic quasilinear quadratic form over $F$ remains anisotropic over $L$ (see \cite[Prop. 5.3]{Hoffmann3}). By the remarks preceding Lemma \ref{LEMbehaviouroverfieldextensions}, it follows that the norm degree is also insensitive to separable extensions. In particular, isotropy indices and norm degrees are insensitive to purely transcendental extensions, something that we shall use without further comment in the sequel. Among inseparable field extensions, a special role is of course played by the inseparable quadratic extensions. If $K$ is a field of characteristic $2$ and $a \in K \setminus K^2$, then we will write $K_a$ for the field $K[X]/(X^2-a)$ ($X$ being an indeterminate). We have here the following:

\begin{lemma} \label{LEMquadraticextensions} Let $\phi$ be an anisotropic quasilinear quadratic form over $F$, and let $a \in F \setminus F^2$. Then:
\begin{enumerate} \item $D(\anispart{(\phi_{F_a})}) = D(\phi_{F_a}) =  D(\phi) + aD(\phi) = D(\pfister{a} \otimes \phi)$;
\item If $\phi_{F_a}$ is anisotropic, then $\pfister{a} \otimes \phi$ is an anisotropic quasilinear quadratic form over $F$ with the same value set as $\anispart{(\phi_{F_a})}$.
\item If $\phi_{F_a}$ is isotropic, then $\witti{0}{\phi_{F_a}} = \mathrm{max}(\mydim{\tau}\;|\; \tau \text{ is a subform of }\phi \text{ divisible by } \pfister{a} )$.  \end{enumerate}
\begin{proof} Statement (1) follows immediately from the additivity property of quasilinear quadratic forms. Now if $\anispart{(\phi_{F_a})}$ is anisotropic, then $D(\anispart{(\phi_{F_a})})$ has dimension $\mydim{\phi}$ as an $F_a^2$-vector space, and hence dimension $2\mydim{\phi}$ as an $F^2$-vector space. Statement (2) then follows from (1). For (3), see \cite[Prop. 5.10]{Hoffmann3}.
\end{proof}
\end{lemma}

\subsection{Function Fields of Quasilinear Quadrics} \label{SUBSECfunctionfields}

If $(V,\phi)$ is a quasilinear quadratic form over $F$ with $\mathrm{lndeg}(\phi) \geq 2$ (e.g., if $\phi$ is anisotropic of dimension $\geq 2$), then the affine hypersurface $X_{\phi}: = \lbrace \phi = 0 \rbrace \subset \mathbb{A}(V)$ is integral, and we write $F(\phi)$ for its function field. Concretely, if $\phi \simeq \form{a} \perp \psi$, and $U$ is the underlying vector space of $\psi$, then $F(\phi)$ is identified with $F(U^\vee)_{a\psi}$. As $X_{\phi}$ is a cone over $X_{\anispart{\phi}}$, $F(\phi)$ is a purely transcendental extension of $F(\anispart{\phi})$. Note that this basic fact, which will be used below, contrasts greatly with the situation for smooth quadrics, where isotropy implies rationality. The form $\phi$ is evidently isotropic after scalar extension to $F(\phi)$. We set $\witti{1}{\phi}: =\witti{0}{\phi_{F(\phi)}} - \witti{0}{\phi}$ and $\phi_1: = \anispart{(\phi_{F(\phi)})}$. In the case where $\phi$ is anisotropic, $\witti{1}{\phi}$ is the smallest non-zero isotropy index attained by $\phi$ over all possible extensions of $F$ (\cite[Cor. 4.8]{Scully1}). If $\phi$ and $\psi$ are anisotropic quasilinear quadratic forms of dimension $\geq 2$ over $F$, then the forms $\phi_{F(\psi)}$ and $\psi_{F(\phi)}$ are simultaneously isotropic precisely when $X_{\phi}$ and $X_{\psi}$ are stably birationally equivalent over $F$ (again, see \cite[\S 2.G]{Scully2}). In this case, we will simply say that $\phi$ and $\psi$ are stably birationally equivalent. For example, if $\phi$ is an anisotropic quasi-Pfister neighbour, then $\phi$ and $\normform{\phi}$ are stably birationally equivalent. In the sequel, a crucial role will be played by the following result on the isotropy of quasilinear quadratic forms over function fields of quadrics: 

\begin{theorem}[{\cite[Thm. 6.4]{Scully1}}] \label{THMkeytool} Let $\phi$ and $\psi$ be anisotropic quasilinear quadratic forms of dimension $\geq 2$ over $F$. If $\phi_{F(\psi)}$ is isotropic, then there exists an anisotropic quasilinear quadratic form $\tau$ of dimension $\witti{0}{\phi_{F(\psi)}}$ over $F(\psi)$ such that $\anispart{(\tau \otimes \psi_1)} \subset \anispart{(\phi_{F(\psi)})}$. 
\end{theorem}

The following corollary records some particular consequences of this result that will also be needed. Several of these were well known prior to \cite{Scully1}. Specifically, part (1) is due to Totaro (\cite[Thm. 5.2]{Totaro}), part (5) is due to Hoffmann (\cite[Thm. 7.7]{Hoffmann3}), and parts (4) and (6) are due to Hoffmann and Laghribi (\cite[Thm. 1.1]{HoffmannLaghribi2}, \cite[Prop. 8.13]{HoffmannLaghribi2}). 

\begin{corollary} \label{CORcorofisotropytheorem} Let $\phi$ and $\psi$ be anisotropic quasilinear quadratic forms of dimension $\geq 2$ over $F$. 
\begin{enumerate} \item If $\phi_{F(\psi)}$ is isotropic, then $\witti{0}{\phi_{F(\psi)}} \leq \mathrm{min}\left(\frac{\mydim{\phi}}{2}, \mydim{\psi} - \witti{1}{\psi} - \mydim{\phi}\right)$. 
\item If $\phi$ and $\psi$ are stably birationally equivalent, then $\psi_1$ is divisible by an $r$-fold quasi-Pfister form, where $r$ is the smallest positive integer for which $\witti{0}{\phi_{F(\psi)}} \leq 2^r$. 
\item $($Values of the invariant $\mathfrak{i}_1)$ $\witti{1}{\psi} \equiv \mydim{\psi} \pmod{2^r}$, where $r$ is the smallest positive integer for which $\witti{1}{\psi} \leq 2^r$.
\item $($Separation theorem$)$ If $\mydim{\phi} \leq 2^s < \mydim{\psi}$ for some integer $s$, then $\phi_{F(\psi)}$ is anisotropic.
\item $\witti{0}{\phi_{F(\psi)}} = \frac{\mydim{\phi}}{2}$ if and only if $\phi$ is divisible by $\normform{\psi}$.
\item If $\phi_{F(\psi)}$ is isotropic, then $\normform{\psi} \subset \normform{\phi}$ and $\mathrm{lndeg}(\phi_{F(\psi)}) = \mathrm{lndeg}(\phi) - 1$.
\end{enumerate}
\begin{proof} (1) Assume $\phi_{F(\psi)}$ is isotropic, and let $\tau$ be as in the statement of Theorem \ref{THMkeytool}. If $a \in D(\psi_1)$, then $D(a\tau) = aD(\tau) \subseteq D(\tau \otimes \psi_1) = D(\anispart{(\tau \otimes \psi_1)})$, and so $a\tau \subset \anispart{(\tau \otimes \psi_1)} \subset \anispart{(\phi_{F(\psi)})}$. Comparing dimensions, we get that $\witti{0}{\phi_{F(\psi)}} \leq \frac{\mydim{\phi}}{2}$. In the same way, we also have that $a\psi_1 \subset \anispart{(\phi_{F(\psi)})}$ for all $a \in D(\tau)$, and dimension comparison then gives that $\witti{0}{\phi_{F(\psi)}} \leq \mydim{\psi} - \witti{1}{\psi} - \mydim{\phi}$. 

(2) If $\phi$ and $\psi$ are stably birationally equivalent, then (1) implies that $\mydim{\psi_1} = \mydim{\psi} - \witti{1}{\psi} = \mydim{\phi} - \witti{1}{\phi} = \mydim{\phi} - \witti{0}{\phi_{F(\psi)}} = \mydim{\anispart{(\phi_{F(\psi)})}}$. The proof of (1) then shows that $a\psi_1 \simeq \anispart{(\phi_{F(\psi)})}$ for all $a \in D(\tau)$. If $b \in D(\phi_{F(\psi)}) \setminus \lbrace 0 \rbrace$, it then follows that $bD(\tau) \subseteq G(\psi_1)$, and so $\simform{(\psi_1)}$ has dimension $\geq 2^r$. 

(3) Applying (2) with $\phi = \psi$, we get that $\psi_1$ is divisible by an $r$-fold quasi-Pfister form, and so $\mydim{\psi} - \witti{1}{\psi} = \mydim{\psi_1} \equiv 0 \pmod{2^r}$.

(4) By (3) we have $\mydim{\psi} - \witti{1}{\psi} \geq 2^s$, and (1) then implies that $\witti{0}{\phi_{F(\psi)}} = 0$.

(5) Since $\psi$ is similar to a subform of $\normform{\psi}$, the latter becomes isotropic over $F(\psi)$. Since an isotropic bilinear Pfister form is split, it follows from Lemma \ref{LEMWittindexrelation} that $\witti{0}{(\normform{\psi})_{F(\psi)}}  \geq \frac{\mydim{\normform{\psi}}}{2}$. In particular, if $\phi$ is divisible by $\normform{\psi}$, then $\witti{0}{\phi_{F(\psi)}} \geq \frac{\mydim{\phi}}{2}$, and so $\witti{0}{\phi_{F(\psi)}} = \frac{\mydim{\phi}}{2}$ by (1). Conversely, suppose that $\witti{0}{\phi_{F(\psi)}} = \frac{\mydim{\phi}}{2}$, and let $\tau$ be as in the proof of (1). If $\witti{0}{\phi_{F(\psi)}} = \frac{\mydim{\phi}}{2}$, then $\mydim{\tau} = \mydim{\anispart{(\phi_{F(\psi)})}}$. The proof of (1) then shows that $a\tau \simeq \anispart{(\phi_{F(\psi)})}$ for all $a \in D(\psi_1)$. As in the proof of (2), this implies that $bD(\psi_1) \subseteq G(\anispart{(\phi_{F(\psi)})})$ for some $b \in F^\times$. Since $G(\anispart{(\phi_{F(\psi)})})$ is a subfield of $F$, we then have that $N(\psi_1) \subseteq G(\anispart{(\phi_{F(\psi)})})$. To prove that $\phi$ is divisible by $\normform{\psi}$, we have to show that $N(\psi) \subseteq G(\phi)$. But if $a \in N(\psi) \setminus \lbrace 0 \rbrace$, then inclusion $N(\psi_1) \subseteq G(\anispart{(\phi_{F(\psi)})})$ implies that the isotropy index of $\anispart{(\pfister{a} \otimes \phi)}$ over $F(\psi)$ is $\mydim{\anispart{(\pfister{a} \otimes \phi)}} - \mydim{\anispart{(\phi_{F(\psi)})}} = \mydim{\anispart{(\pfister{a} \otimes \phi)}} - \frac{\mydim{\phi}}{2}$. By (1) (with $\phi$ replaced by $\anispart{(\pfister{a} \otimes \phi)}$), it then follows that $\mydim{\anispart{(\pfister{a} \otimes \phi)}} = \mydim{\phi}$, which means that $a\phi \simeq \phi$, i.e., $a \in G(\phi)$.

(6) As noted in the proof of (5), the isotropy of $\phi_{F(\psi)}$ implies that $\witti{0}{(\normform{\phi})_{F(\psi)}} = \frac{\mydim{\normform{\phi}}}{2}$. Applying (5), we then get that $\normform{\psi} \subset \normform{\phi}$. Moreover, the equality $\mathrm{lndeg}(\phi_{F(\psi)}) = \mathrm{lndeg}(\phi) - 1$ also holds since $\normform{(\phi_{F(\psi)})} \simeq \anispart{((\normform{\phi})_{F(\psi)})}$.  \end{proof} \end{corollary}

\begin{remarks} \label{REMSisotropytheoremcorollaries} \begin{enumerate}[leftmargin=*] \item Part (3) is the quasilinear version of an important theorem of Karpenko and Primozic describing the possible values of the invariant $\mathfrak{i}_1$ for non-singular quadratic forms (see \cite{Karpenko1, Primozic}). If $s$ is the unique integer for which $2^s < \mydim{\psi} \leq 2^{s+1}$, then it implies that $\witti{1}{\psi} \leq \mydim{\psi} - 2^s$. By parts (2) and (5), equality holds here in the case where $\psi$ is a quasi-Pfister neighbour (though the converse is not true in general).
\item Applying (6) with $\phi = \psi$, we also get that $\mathrm{lndeg}(\psi_1) = \mathrm{lndeg}(\psi) -1$, a fact that will be used several times below. 
\end{enumerate}
\end{remarks}

By Lemma \ref{LEMWittindexrelation}, the hypothesis in part (5) of the corollary is satisfied when $\phi$ is the quadratic form associated to a anisotropic symmetric bilinear form that splits over $F(\psi)$. One can then use (5) to deduce the following important result due to Laghribi that will also be needed:

\begin{theorem}[{\cite[Thm. 1.2]{Laghribi1}}] \label{THMWittkernel} Let $\mathfrak{b}$ be an anisotropic symmetric bilinear form over $F$, and let $\psi$ be an anisotropic quasilinear quadratic form over $F$. If $\mathfrak{b}_{F(\psi)}$ is split, then there exist scalars $x_1,\hdots,x_r \in F^\times$ and bilinear Pfister forms $\mathfrak{b}_1,\hdots,\mathfrak{b}_r$ over $F$ such that $\mathfrak{b} = x_1\mathfrak{b}_1 \perp \cdots \perp x_r\mathfrak{b}_r$ and $\phi_{\mathfrak{b}_i} \simeq \normform{\psi}$ for all $1 \leq i \leq r$. 
\end{theorem}

Note that the Pfister forms $\mathfrak{b}_1,\hdots,\mathfrak{b}_r$ here need not be isometric here, in contrast to the situation over fields of characteristic not 2. We now make some simple observations on bilinear lifts of quasilinear quadratic forms that will be used later. 

\section{Bilinear Lifts of Quasilinear Quadratic Forms} \label{SECbilinearlifts}

\subsection{Lifts of Anisotropic Quasi-Pfister Forms} As indicated above, an anisotropic quasilinear quadratic form will generally lift to infinitely many non-isometric symmetric bilinear forms. Nevertheless, Theorem \ref{THMWittkernel} calls for an observation regrading the bilinear Pfisters forms that lift a given anisotropic quasi-Pfister form. Following \cite{Laghribi1}, let us write $\mathscr{A}(\pi)$ for the set of consisting of the bilinear Pfister forms over $F$ whose associated quadratic form is isometric to a given quasi-Pfister form $\pi$. When $\pi$ is anisotropic, the following lemma says that the classes of the elements of $\mathscr{A}(\pi)$, together with the zero element, form an additive subgroup of $W(F)$:

\begin{lemma} \label{LEMliftsofquasiPfisters} Let $\pi$ be an anisotropic quasi-Pfister form over $F$, and let $\mathfrak{b}, \mathfrak{c} \in \mathscr{A}(\pi)$. If $\mathfrak{b} \not \simeq \mathfrak{c}$, then $\anispart{(\mathfrak{b} \perp \mathfrak{c})} \in \mathscr{A}(\pi)$. 
\begin{proof} Let $n$ be the foldness of $\pi$, and consider the form $\mathfrak{d} := \anispart{(\mathfrak{b} \perp \mathfrak{c})}$. Since $\mathfrak{b}, \mathfrak{c} \in \mathscr{A}(\pi)$, we have $\phi_{\mathfrak{d}} \subset \anispart{(\pi \perp \pi)} \simeq \pi$, and so $\mydim{\mathfrak{d}} \leq \mydim{\pi} = 2^n$. Note, however that $[\mathfrak{d}] \in I^n(F)$. In particular, if the former dimension inequality is strict, then the Hauptsatz (see \S \ref{SECintroduction} or Theorem \ref{THMHauptsatz} below) tells us that $\mathfrak{d}$ is trivial, and so $\mathfrak{b} \simeq \mathfrak{c}$. If not, then the same result tells us that $\mathfrak{d}$ is a general bilinear Pfister form with $\phi_{\mathfrak{d}} \simeq \pi$. But a general bilinear Pfister form that represents $1$ is Pfister, so we in fact have that $\mathfrak{d} \in \mathscr{A}(\pi)$ in this case. \end{proof} \end{lemma}

\subsection{Lifts Inside a Fixed Form and a Partial Splitting Construction} In general, the non-uniqueness of bilinear lifts of anisotropic quasilinear quadratic forms may be resolved by asking that for a lift that is embedded into a fixed anisotropic form:

\begin{lemma} \label{LEMliftingsubforms} Let $\mathfrak{b}$ be an anisotropic symmetric bilinear form over $F$, and let $\phi$ be  a subform of the quasilinear quadratic form $\phi_{\mathfrak{b}}$. Up to isometry, there then exists a unique subform $\mathfrak{c}$ of $\mathfrak{b}$ such that $\phi_{\mathfrak{c}} \simeq \phi$. 
\begin{proof} Let $V$ be the underlying vector space of $\mathfrak{b}$. If $\phi$ is isometric to the restriction of $\phi_{\mathfrak{b}}$ to a subspace $U$ of $V$, then $\phi_{\mathfrak{b}|_U} \simeq \phi$. For uniqueness, suppose that $\mathfrak{c}$ and $\mathfrak{c'}$ are subforms of $\mathfrak{b}$ with $\phi_{\mathfrak{c}} \simeq \phi_{\mathfrak{c'}}$. Without loss of generality, we can assume that $\mathfrak{c} = \mathfrak{b}|_U$ and $\mathfrak{c}' = \mathfrak{b}|_{U'}$ for subspaces $U$ and $U'$ of $V$ with $\mathrm{dim}(U) = \mathrm{dim}(U')$. Suppose that $U \neq U'$. Then there exists a vector $u \in U \setminus U'$. Since $\phi_{\mathfrak{c}} \simeq \phi_{\mathfrak{c'}}$, we can find a vector $u' \in U'$ such that $\phi_{\mathfrak{c}'}(u') = \phi_{\mathfrak{c}}(u)$. But then $\phi_{\mathfrak{b}}(u + u') = \phi_{\mathfrak{c}}(u) + \phi_{\mathfrak{c}'}(u') = 0$, and since $u + u' \neq 0$, this contradicts the anisotropy of $\mathfrak{b}$. Thus, $U = U'$, and so the uniqueness claim holds. \end{proof}
\end{lemma}

In the sequel, we will need to construct function fields of quasilinear quadrics that achieve a desired degree of splitting in a given quasilinear quadratic form. In \cite[\S 4.a]{Scully2}, it was noted that the existence part of Lemma \ref{LEMliftingsubforms} sometimes provides an effective means of doing this. The basic point here is the following simple observation:

\begin{lemma}[{see \cite[Cor. 4.12]{Scully2}}] \label{LEMisotropylemma} Let $\mathfrak{b}$, $\mathfrak{c}$ and $\mathfrak{d}$ be anisotropic symmetric bilinear forms over $F$ with $\mathfrak{b} \simeq \mathfrak{c} \perp \mathfrak{d}$. If $\psi$ is an anisotropic quasilinear quadratic form of dimension $\geq 2$ over $F$ such that $\mathfrak{b}_{F(\psi)}$ is split, then  $\mathfrak{i}_{\mathbb{H}}(\mathfrak{b}_{F(\psi)}) = 0$ and $\witti{0}{(\phi_{\mathfrak{c}})_{F(\psi)}} = \witti{0}{(\phi_{\mathfrak{d}})_{F(\psi)}} + \frac{\mydim{\mathfrak{c}} - \mydim{\mathfrak{d}}}{2}$.
\begin{proof} Since $\mathfrak{b}_{F(\psi)}$ is split, we have $[\mathfrak{c}_{F(\psi)}]=[\mathfrak{d}_{F(\psi)}]$ in $W(F(\psi))$, and so $\windex{\mathfrak{c}_{F(\psi)}} = \windex{\mathfrak{d}_{F(\psi)}} + \frac{\mydim{\mathfrak{c}} - \mydim{\mathfrak{d}}}{2}$. By Lemma \ref{LEMWittindexrelation}, proving the second equality in the statement then amounts to showing that $\mathfrak{i}_{\mathbb{H}}(\mathfrak{c}_{F(\psi)}) = \mathfrak{i}_{\mathbb{H}}(\mathfrak{d}_{F(\psi)}) = 0$. But the invariant $\mathfrak{i}_{\mathbb{H}}$ is clearly subadditive with respect to orthogonal sums, so if we can show that $\mathfrak{i}_{\mathbb{H}}(\mathfrak{b}_{F(\psi)}) = 0$, then the second equality in the statement will also folllow. By Lemma \ref{LEMWittindexrelation}, we have
$$ \mathfrak{i}_{\mathbb{H}}(\mathfrak{b}_{F(\psi)}) = \witti{0}{(\phi_{\mathfrak{b}})_{F(\psi)}} - \windex{\mathfrak{b}_{F(\psi)}} = \witti{0}{(\phi_{\mathfrak{b}})_{F(\psi)}} - \frac{\mydim{\mathfrak{b}}}{2} = \witti{0}{(\phi_{\mathfrak{b}})_{F(\psi)}} - \frac{\mydim{\phi_{\mathfrak{b}}}}{2}. $$
Since the isotropy index attained by an anisotropic quasilinear quadratic form after extension to the function field of a quasilinear quadric cannot exceed half its dimension (see part (1) of Corollary \ref{CORcorofisotropytheorem}), it then follows that $\mathfrak{i}_{\mathbb{H}}(\mathfrak{b}_{F(\psi)}) = 0$, as desired. \end{proof}
\end{lemma}

We will be interested in the following application:

\begin{proposition} \label{PROPcomplementarytoneighbours} Let $n$ be a non-negative integer, $\eta$ an anisotropic $(n+1)$-fold quasi-Pfister form over $F$, and $\phi$ a subform of $\eta$ of codimension at least $2$. Then there exists a quasilinear quadratic form $\psi$ over $F$ such that $\eta \simeq \phi \perp \psi$ and $\witti{0}{\phi_{F(\psi)}} = \witti{1}{\psi} + \mydim{\phi} - 2^n$. 
\begin{proof} Choose a bilinear Pfister form $\mathfrak{b} \in \mathscr{A}(\eta)$. By Lemma \ref{LEMliftingsubforms}, there exists an orthogonal decomposition $\mathfrak{b} = \mathfrak{c} \perp \mathfrak{d}$ with $\phi_{\mathfrak{c}} = \phi$. Set $\psi: = \phi_{\mathfrak{d}}$. We then clearly have that $\eta \simeq \phi \perp \psi$. Moreover, the form $\mathfrak{b}$ becomes isotropic over $F(\psi)$, and hence splits over $F(\psi)$ on account of being Pfister. By Lemma \ref{LEMisotropylemma}, it then follows that
$$ \witti{0}{\phi_{F(\psi)}} = \witti{1}{\psi} + \frac{\mydim{\mathfrak{c}} - \mydim{\mathfrak{d}}}{2} = \witti{1}{\psi} + \frac{2\mydim{\phi} - 2^{n+1}}{2} = \witti{1}{\psi} + \mydim{\phi} - 2^n, $$
and so $\psi$ has the desired properties.  \end{proof}
\end{proposition}

Note here that any given subform of a quasilinear quadratic form will generally have infinitey many non-isometric complementary forms, so $\psi$ must be chosen carefully in the situation of the proposition. The statement is specifically of interest in the case where $\mydim{\phi}> 2^n$, meaning that $\phi$ is a quasi-Pfister neighbour with $\normform{\phi} =\eta$. Here, the proposition leads to a short proof of Theorem \ref{THMdimensiongaptheorem} in the characteristic-2 case, as we now explain. 

\section{Dimension Gaps in Characteristic 2} \label{SECgaps}

\subsection{Elementary Proof of the Dimension Gap Theorem in Characteristic 2} Our goal here is to give a short and elementary proof of Theorem \ref{THMdimensiongaptheorem} in the characteristic-2 case. For this, we will need the Hauptsatz of Arason and Pfister discussed in \S \ref{SECintroduction}. In the interest of completeness, we quickly explain how this follows from Theorem \ref{THMWittkernel} above, which will be a basic tool in all that follows. The proof further determines the anisotropic forms of dimension $2^n$ representing an element of $I^n(F)$, which we will also need.

\begin{theorem}[Hauptsatz] \label{THMHauptsatz} Let $n$ be a positive integer, and let $\mathfrak{b}$ be a non-zero anisotropic symmetric bilinear form over $F$ that represents an element of $I^n(F)$. Then $\mydim{\mathfrak{b}} \geq 2^n$. Moreover, if equality holds, then $\mathfrak{b}$ is a general $n$-fold bilinear Pfister form.
\begin{proof} For the first statement, we argue by induction on the smallest integer $r$ with the property that $[\mathfrak{b}]$ can be written as the sum of the classes of $r$ general $n$-fold bilinear Pfister forms in $W(F)$. If $r =1$, there is nothing to show. Suppose now that $r \geq 2$, and let $\mathfrak{b}_1,\hdots,\mathfrak{b}_r$ be general $n$-fold bilinear Pfister forms over $F$ such that $[\mathfrak{b}] = [\mathfrak{b}_1] + \cdots + [\mathfrak{b}_r]$ in $W(F)$. Set $\psi: = \phi_{\mathfrak{b}_r}$. Then $\mathfrak{b}_r$ splits over $F(\psi)$, and so $[\anispart{(\mathfrak{b}_{F(\psi)})}]$ is expressible as the sum of less than $r$ general $n$-fold bilinear Pfister forms in $W(F(\psi))$. Applying the induction hypothesis (over $F(\psi
)$), we get that $\mydim{\anispart{(\mathfrak{b}_{F(\psi)})}}$ is either $0$ or $\geq 2^n$. In the second case, we clearly have that $\mydim{\mathfrak{b}} \geq 2^n$. In the first case, however, the same conclusion is forced by Theorem \ref{THMWittkernel} (which tells us that $\mathfrak{b}$ decomposes as an orthogonal sum of general $n$-fold bilinear Pfister forms). The first statement therefore holds. At the same time, if $\mydim{\mathfrak{b}} = 2^n$, then the first statement implies that $\mathfrak{b}$ splits over $F(\phi)$, where $\phi := \phi_{\mathfrak{b}}$. By another application of Theorem \ref{THMWittkernel}, it follows that $\mathfrak{b}$ decomposes as an orthogonal sum of general bilinear Pfister forms of dimension $\mydim{\normform{\phi}} = 2^{\mathrm{lndeg}(\phi)}$. Since $\mydim{\phi} = 2^n$, we have $\mathrm{lndeg}(\phi) \geq n$, and so $\mathfrak{b}$ is in fact a general $n$-fold bilinear Pfister form.
\end{proof} \end{theorem}

We are now ready to prove the dimension gap theorem in characteristic 2:

\begin{theorem}[Dimension Gap Theorem] \label{THMgaptheoreminchar2} Let $n$ be a positive integer, and let $\mathfrak{b}$ be an anisotropic symmetric bilinear form of dimension $<2^{n+1}$ over $F$. If $\mathfrak{b}$ represents an element of the ideal $I^n(F)$, then $\mydim{\mathfrak{b}} = 2^{n+1} - 2^i$ for some integer $1 \leq i \leq n+1$. 
\begin{proof} By the Hauptsatz, we can assume that $\mydim{\mathfrak{b}} > 2^n$, say $\mathfrak{b} = 2^{n+1} - m$ with $1 \leq m < 2^n$. Set $\phi: = \phi_{\mathfrak{b}}$. Since $\mydim{\phi} = 2^{n+1} - m > 2^n$, we have $\mathrm{lndeg}(\phi) \geq n+1$, say $\mathrm{lndeg}(\phi) = n+1 + s$ for a non-negative integer $s$. Suppose that $s \geq 1$. The form $\normform{\phi}$ then has dimension $\geq 2^{n+2}$. Since $\mydim{\phi} = 2^{n+1}-m  < 2^{n+1}$, the separation theorem (Corollary \ref{CORcorofisotropytheorem} (4)) tells us that $\phi$ remains anisotropic over $F(\normform{\phi})$. At the same time, Corollary \ref{CORcorofisotropytheorem} (6) and Remark \ref{REMSisotropytheoremcorollaries} (2) give that $\mathrm{lndeg}(\phi_{F(\normform{\phi})}) = \mathrm{lndeg}((\normform{\phi})_{F(\normform{\phi})}) = \mathrm{lndeg}(\normform{\phi})-1 = \mathrm{lndeg}(\phi) - 1$. Thus, by passing from $F$ to $F(\normform{\phi})$, we can lower the value of $\mathrm{lndeg}(\phi)$ by $1$ without making $\phi$ isotropic. Repeating this procedure a further $s-1$ times, we find an extension $L$ of $F$ such that $\phi_L$ is anisotropic and $\mathrm{lndeg}(\phi_L) = n+1$. To prove the theorem, we can therefore assume that $s = 0$, i.e., that $\phi$ is a quasi-Pfister neighbour. By Proposition \ref{PROPcomplementarytoneighbours} (applied with $\eta = \normform{\phi}$), we can then find a quasilinear quadratic form $\psi$ over $F$ such that $\mydim{\psi} = m$ and $\witti{0}{\phi_{F(\psi)}} = \witti{1}{\psi} + 2^n - m > 2^n - m$. By Lemma \ref{LEMWittindexrelation}, it follows that $2\windex{\mathfrak{b}_{F(\psi)}} \geq \windex{\mathfrak{b}_{F(\psi)}} + \mathfrak{i}_{\mathbb{H}}(\mathfrak{b}_{F(\psi)}) = \witti{0}{\phi_{F(\psi)}} > 2^n - m$, and so $\mydim{\anispart{(\mathfrak{b}_{F(\psi)})}} = \mydim{\mathfrak{b}} - 2\windex{\mathfrak{b}_{F(\psi)}} < (2^{n+1}-m)-(2^n - m) = 2^n$. Since $[\anispart{(\mathfrak{b}_{F(\psi)})}] \in I^n(F(\psi))$, the Hauptsatz then gives that $\anispart{(\mathfrak{b}_{F(\psi)})}$ is zero, i.e., that $\mathfrak{b}$ splits over $F(\psi)$. By Theorem \ref{THMWittkernel}, it then follows that $2^{n+1} - m = \mydim{\mathfrak{b}}$ is divisible by $2^{\mathrm{lndeg}(\psi)}$. Since $\mydim{\psi} = m$, this gives that $2^{n+1} -m \equiv 0 \pmod{2^i}$, where $i$ is the least positive integer with $m \leq 2^i$. Since $m < 2^n$, we must then have that $m = 2^i$, and so we're done.
\end{proof}
\end{theorem}

The proof given here exposes a fundamental difference between the characteristic-not-2 and characteristic-2 cases:

\begin{remarks} \label{REMSnonrigidity} Let $K$ be a field of any characteristic, let $n$ be a positive integer, and let $1 \leq i \leq n-1$. Suppose that $\mathfrak{b}$ is anisotropic symmetric bilinear form of dimension $2^{n+1} - 2^i$ that represents an element of $I^n(K)$ (such forms necessarily exist -- see \S \ref{SUBSECoptimality} below). Let us write $\witti{1}{\mathfrak{b}}$ for the Witt index of $\mathfrak{b}$ after scalar extension to the field $K(\phi)$, where $\phi$ is the quadratic form associated to $\mathfrak{b}$. The dimension of the anisotropic part of $\mathfrak{b}$ over $K(\phi)$ is $\mydim{\mathfrak{b}} - 2\witti{1}{\mathfrak{b}}$, and so Theorem \ref{THMdimensiongaptheorem} gives that $\witti{1}{\mathfrak{b}} = 2^{j-1} - 2^{i-1}$ for some $i \leq j \leq n+1$. In the characteristic-not-$2$ setting, the equality $j=i$ is forced by Karpenko's theorem on the values of $\mathfrak{i}_1$ for non-singular quadratic forms (see Remark \ref{REMSisotropytheoremcorollaries} (1)). In other words, if we run over the Knebusch splitting tower of $\mathfrak{b}$ in this setting (see \cite[\S 25]{EKM}), then the dimensions of the anisotropic parts of $\mathfrak{b}$ that manifest are precisely the integers $2^{n+1} - 2^j$ with $i \leq j \leq n$. The splitting pattern of $\mathfrak{b}$ is therefore completely determined, and it is also known that various other discrete (motivic) invariants of $\mathfrak{b}$ behave in a similarly rigid way in this setting (see, e.g., \cite[\S 82]{EKM}). In the case where $\mathrm{char}(K) = 2$, however, the situation is different. Indeed, although Corollary \ref{CORcorofisotropytheorem} (5) provides the direct analogue of Karpenko's theorem on $\mathfrak{i}_1$ for the quasilinear quadratic form $\phi$, this now only constrains $\mathfrak{i}_1(\mathfrak{b})$ up to knowing the integer $\mathfrak{i}_{\mathbb{H}}(\mathfrak{b}_{K(\phi)})$ (Lemma \ref{LEMWittindexrelation}). If $\mathfrak{i}_{\mathbb{H}}(\mathfrak{b}_{K(\phi)}) = 0$, we again get that $j=i$. It is in general possible, however, that $\mathfrak{i}_{\mathbb{H}}(\mathfrak{b}_{K(\phi)}) \neq 0$, forcing $j$ to be greater than $i$. In fact, we have seen in our proof of Theorem \ref{THMgaptheoreminchar2} that the form $\phi$ may be a quasi-Pfister neighbour. In this case, we have $\witti{1}{\phi} = 2^n - 2^i$, and so $\witti{1}{\mathfrak{b}} \geq 2^{n-1} - 2^{i-1}$. In other words, the integer $j$ may be as large as $n$ here. Since $\witti{1}{\phi}$ is the minimal non-zero isotropy index that $\phi$ can attain over any extension of $K$, it follows that the splitting behaviour of the low-dimensional forms representing elements of $I^n(K)$ is not determined by dimension, and is thus more complicated here than in the characteristic-not-$2$ setting. \end{remarks}

\subsection{Optimality} \label{SUBSECoptimality} It is easy to see that all dimensions permitted in the statement of Theorem \ref{THMgaptheoreminchar2} are actually realizable over a suitable field. For a bilinear Pfister form $\mathfrak{b}$ of dimension $\geq 2$ over a field of characteristic $2$, let us write $\mathfrak{b}'$ for its pure part, i.e., the complementary subform of $\form{1}_b$ in $\mathfrak{b}$. We then have the following:

\begin{lemma} \label{LEMoptimalityinlowdimension} Let $i \leq n$ be positive integers, let $X_1,\hdots,X_{i-1}, Y_i,\hdots,Y_n,Z_i,\hdots,Z_n$ be indeterminates, and set $K: = F(X_1,\hdots,X_{i-1},Y_i,\hdots,Y_n,Z_1,\hdots,Z_n)$. Then
$$ \mathfrak{b}: = \pfister{X_1,\hdots,X_{i-1}}_b \otimes \left(\pfister{Y_i,\hdots,Y_n}_b' \perp \pfister{Z_i,\hdots,Z_n}_b'\right) $$
is an anisotropic symmetric bilinear form of dimension $2^{n+1} - 2^i$ over $K$ with $[\mathfrak{b}] \in I^n(K)$. 
\begin{proof} We have $[\mathfrak{b}] = [\pfister{X_1,\hdots,X_{i-1},Y_i,\hdots,Y_n}_b] +  [\pfister{X_1,\hdots,X_{i-1},Z_i,\hdots,Z_n}_b]$ in $W(K)$, and so $[\mathfrak{b}] \in I^n(K)$. It then only remains to check that $\mathfrak{b}$ is anisotropic, but this follows from the fact that $\lbrace X_1,\hdots,X_{i-1},Y_i,\hdots,Y_n,Z_i,\hdots,Z_n \rbrace$ is a $2$-independent subset of $K$. \end{proof}
\end{lemma}

As far as the optimality of Theorem \ref{THMgaptheoreminchar2} goes, the interesting question is therefore whether further dimension gaps can appear if we look at the anisotropic symmetric bilinear forms of even dimension $\geq 2^{n+1}$ that represent an element of $I^n(F)$. In the characteristic-not-2 setting, Vishik has shown that no further gaps can appear using motivic methods (see \cite[Prop. 82.7]{EKM}). In characteristic 2, however, the situation seems to be less clear. Given a symmetric bilinear form $\mathfrak{b}$ over a field $K$, let us write $\mathrm{Dim}(\mathfrak{b})$ for the set $\lbrace \mydim{\anispart{(\mathfrak{b}_L)}}\;|\; L/K \text{ a field extension} \rbrace$. Let $m$ and $n$ be positive integers with $m \geq 3$, let $k$ be a field (of any characteristic), and consider the form
$$\mathfrak{c}: = X_1 \pfister{X_{11},\hdots,X_{1n}}_b \perp \cdots \perp X_m \pfister{X_{m1},\hdots,X_{mn}}_b, $$
of dimension $m2^n$ over the rational function field $K: = k(X_i,X_{ij}\;|\; 1 \leq i \leq m, 1 \leq j \leq n)$. Clearly $\mathfrak{c}$ represents an element of $I^n(K)$, and every symmetric bilinear form expressible as a sum of $m$ general $n$-fold bilinear Pfister forms over an extension of $k$ is a specialization of $\mathfrak{c}$. With a view to the optimality problem, it is therefore natural to try to determine the set $\mathrm{Dim}(\mathfrak{c})$. When $\mathrm{char}(k) \neq 2$, Vishik showed specifically that $\mathrm{Dim}(\mathfrak{c})$ contains the set $\lbrace 2r\;|\; 2^n \leq r \leq m2^{n-1} \rbrace$, thereby settling the problem positively. The arguments here rely on the fact that in order to construct the set $\mathrm{Dim}(\mathfrak{c})$ in this setting, one only needs to consider fields in the Knebusch splitting tower of $\mathfrak{c}$. When $\mathrm{char}(k) =2$, however, the Knebusch splitting tower need not be universal among the splitting towers of $\mathfrak{c}$, and may therefore not realize all values of the set $\mathrm{Dim}(\mathfrak{c})$. Moreover, as far as our problem goes, it is in fact impossible for the component of $\mathrm{Dim}(\mathfrak{c})$ determined by the Knebusch splitting tower to include all elements of the set $\lbrace 2r\;|\; 2^n \leq r \leq m2^{n-1} \rbrace$. Indeed, as we move from one step of the tower to the next, the norm degree of the quasilinear quadratic form $\phi_{\mathfrak{c}}$ drops by $1$. As the initial norm degree is only $mn + m-1$, it will eventually become small enough that a non-trivial jump in the anisotropic dimension of $\mathfrak{c}$ will be forced (see Remark \ref{REMSisotropytheoremcorollaries} (1)). A direct adaptation of Vishik's construction is therefore not feasible here. Nevertheless, it is still possible to use splitting towers to obtain some partial restrictions. For the remainder of this section, let us set
$$ \mathrm{Dim}(I^n,F): = \lbrace \mydim{\anispart{\mathfrak{b}}}\;|\; [\mathfrak{b}] \in I^n(K) \text{ for some extension }K/F \rbrace. $$
We have the following simple lemma:

\begin{lemma} \label{LEMadding2^n} Let $n$ be a positive integer. If $d \in \mathrm{Dim}(I^n,F)$, then $d+2^n \in \mathrm{Dim}(I^n,F)$. 
\begin{proof} Since $d \in \mathrm{Dim}(I^n,F)$, there exists an extension $K$ of $F$ and an anisotropic symmetric bilinear form $\mathfrak{b}$ of dimension $d$ over $K$ such that $[\mathfrak{b}] \in I^n(K)$. Set $L: = K(X_0,\hdots,X_n)$, where $X_0,\hdots,X_n$ are indeterminates, and consider the form $\mathfrak{c} : = \mathfrak{b} \perp X_0 \pfister{X_1,\hdots,X_n}_b$ of dimension $d + 2^n$ over $L$. It is then clear that $[\mathfrak{c}] \in I^n(L)$, so to prove the lemma, it suffices to show that $\mathfrak{c}$ is anisotropic. But $\pfister{X_1,\hdots,X_n}_b$ is anisotropic since $\lbrace X_1,\hdots,X_n \rbrace$ is a $2$-independent subset of $L$, and the anisotropy of $\mathfrak{c}$ then follows from \cite[Cor. 19.6]{EKM}. \end{proof}
\end{lemma}

The basic issue is therefore to determine which even integers in the interval $[2^{n+1},2^{n+1} + 2^n)$ lie in $\mathrm{Dim}(I^n,F)$. We can show here the following:

\begin{lemma} \label{LEMtechnicaloptimalitylemma} Let $n$ be a positive integer, let $1 \leq i \leq n-1$, and let $0 \leq a, b \leq n-i+1$ be such that $a +b < n$. Then there exist a field extension $K$ of $F$ and an anisotropic symmetric bilinear form $\mathfrak{b}$ of dimension $2^{n+1} + 2(2^{n-1} - 2^{i-1} - 2^a - 2^b + 2)$ over $K$ such that $[\mathfrak{b}] \in I^n(K)$. Moreover, when $n \geq 5$ and $i = n-1$, we can choose $\mathfrak{b}$ so that $\witti{1}{\phi_{\mathfrak{b}}}= 1$. 
\begin{proof} Set $K: = F(X_1,\hdots,X_{i-1},Y_i,\hdots,Y_n,Z_i,\hdots,Z_n,W_{a+b+1},\hdots,W_n)$, where the $X_j$, $Y_j$, $Z_j$ and $W_j$ are indeterminates. Set
$$ \mathfrak{c} : = \pfister{X_1,\hdots,X_{i-1}}_b \otimes \left(\pfister{Y_i,\hdots,Y_n}_b' \perp \pfister{Z_i,\hdots,Z_n}_b'\right), $$
$$ \mathfrak{d} : = \pfister{Y_i,\hdots, Y_{i+a-1},Z_{i},\hdots,Z_{i+b-1},W_{a+b+1},\hdots,W_n}_b, $$
and $ \mathfrak{b} : = \anispart{(\mathfrak{c} \perp \mathfrak{d})}$. By Lemma \ref{LEMoptimalityinlowdimension}, $\mathfrak{c}$ is an anisotropic form of dimension $2^{n+1} - 2^i$ with $[\mathfrak{c}] \in I^n(K)$. Since $\mathfrak{d}$ is an $n$-fold Pfister form, we then also have that $[\mathfrak{d}] \in I^n(K)$. Using that $\lbrace X_1,\hdots,W_n \rbrace$ is a $2$-independent subset of $K$, one also readily checks that $\mydim{\mathfrak{b}} = 2^{n+1} - 2^i + 2^n - 2x$, where $x:= 2^a + 2^b -2$. In fact, $\pfister{Y_i,\hdots,Y_{i+a-1}}' \perp \pfister{Z_i,\hdots,Z_{i+a-1}}'$ is a common $x$-dimensional subform of $\mathfrak{c}$ and $\mathfrak{d}$, and the orthogonal sum of the complementary subforms is easily seen to be anisotropic, and hence equal to $\mathfrak{b}$. This proves the first statement. Suppose now that $n \geq 5$ and $i = n-1$. Since $i = n-1$, we then have that $a,b \leq 2$, and so $\mydim{\mathfrak{b}} = 2^{n+1} +2^{n-1} -2x$ with $x \leq 6$. Since $n \geq 5$, we in particular have that $2^{n+1} < \mydim{\mathfrak{b}} \leq 2^{n+1} + 2^{n-1}$. By construction, we also have that $\mathrm{lndeg}(\phi_{\mathfrak{b}}) = 2n + 2 - (a+b) \geq n+3$, and so $\phi_{\mathfrak{b}}$ is not a quasi-Pfister neighbour. For ease of notation, let us now set $\phi: = \phi_{\mathfrak{b}}$. We wish to show that $\witti{1}{\phi} = 1$. Suppose, for the sake of contradiction, that this does not hold. By Corollary \ref{CORcorofisotropytheorem}, the quadratic form $\phi_1$ is then divisible by a binary form, say $\pfister{a}$ with $a \in K(\phi)\setminus K(\phi)^2$. We claim that $a \in D(\pfister{X_1,\hdots,X_{n-2}}_{K(\phi)})$. Let us start by observing that $(\phi_{\mathfrak{c}})_{K(\phi)}$ is an anisotropic non-quasi-Pfister neighbour. Indeed, $\mathrm{lndeg}(\phi_{\mathfrak{c}}) = n+2$ by construction. Since $\mathrm{lndeg}(\phi) \geq n+3$, it follows from part (6) of Corollary \ref{CORcorofisotropytheorem} that $\normform{(\phi_{\mathfrak{c}})}$ remains anisotropic over $K(\phi)$. Since $\phi_{\mathfrak{c}}$ is similar to a subform of $\normform{(\phi_{\mathfrak{c}})}$, it must then also remain anisotropic over $K(\phi)$. Furthermore, we have $\mathrm{lndeg}((\phi_{\mathfrak{c}})_{K(\phi)}) = \mathrm{lndeg}(\phi_{\mathfrak{c}}) = n+2$, and since $\mydim{\mathfrak{c}} = 2^n + 2^{n-1} <2^{n+1}$, $(\phi_{\mathfrak{c}})_{K(\phi)}$ is therefore not a quasi-Pfister neighbour. As the next step, let us now show that $\mathfrak{c}$ splits over $K(\phi)_a$. Suppose, to the contrary, that this is not the case. Since $\phi_1$ is divisible by $\pfister{a}$, its isotropy index over $F(\phi)_a$ is half its dimension (see part (5) of Corollary \ref{CORcorofisotropytheorem}), and so
\begin{equation} \label{eq0} \mydim{\anispart{(\phi_{K(\phi)_a})}} = 2^n + 2^{n-2} - x - \frac{\witti{1}{\phi}}{2}. \end{equation}
Recall, however that $\mathfrak{b}$ contains a codimension-$x$ subform of $\mathfrak{c}$, and so this shows that $\mydim{\anispart{((\phi_{\mathfrak{c}})_{K(\phi)_a})}} < 2^n + 2^{n-2}$. Since $\mathfrak{c}$ is not split over $K(\phi)_a$, however, the quasi-Pfister form $\pfister{X_1,\hdots,X_{n-2}}$ remains anisotropic over $K(\phi)_a$, and is hence a divisor of $\anispart{((\phi_{\mathfrak{c}})_{K(\phi)_a})}$ by Lemma \ref{LEMbehaviouroverfieldextensions} (2). The above remarks therefore give that $\mydim{\anispart{((\phi_{\mathfrak{c}})_{K(\phi)_a})}} \leq 2^n$. At the same time, $\anispart{((\phi_{\mathfrak{c}})_{K(\phi)_a})}$ contains the quadratic form associated to $\anispart{(\mathfrak{c}_{K(\phi)_a})}$ as a subform. Since $\mathfrak{c}_{K(\phi)_a}$ is not split, it follows from Theorem \ref{THMHauptsatz} that this subform is similar to an $n$-fold quasi-Pfister form. But Corollary \ref{CORcorofisotropytheorem} (6) then gives that $\mathrm{lndeg}((\phi_{\mathfrak{c}})_{K(\phi)}) = \mathrm{lndeg}((\phi_{\mathfrak{c}})_{K(\phi)_a}) + 1 = n + 1$, contradicting our earlier observation that $(\phi_{\mathfrak{c}})_{K(\phi)}$ is not a quasi-Pfister neighbour. We can therefore conclude that $\mathfrak{c}_{K(\phi)_a}$ is indeed split. Since $\mathfrak{c}_{K(\phi)}$ is anisotropic, Theorem \ref{THMWittkernel} then implies that $(\phi_{\mathfrak{c}})_{K(\phi)}$ is divisible by $\pfister{a}$. In other words, we have $a \in G((\phi_{\mathfrak{c}})_{K(\phi)})$. Now, if an anisotropic quasilinear quadratic form of dimension $2^n + 2^{n-1}$ is divisible by an $(n-1)$-fold quasi-Pfister form, then it is has norm degree $n+1$ and is hence a quasi-Pfister neighbour. Since $(\phi_{\mathfrak{c}})_{K(\phi)}$ is not a quasi-Pfister neighbour, it follows that $\pfister{X_1,\hdots,X_{n-2}} \simeq \simform{((\phi_{\mathfrak{c}})_{K(\phi)})}$, and so $a \in D(\pfister{X_1,\hdots,X_{n-2}}_{K(\phi)})$, as claimed. Now, since $\mydim{\phi} > 2^{n+1}$, we have $\mydim{\phi_1} \geq 2^{n+1}$ by Remark \ref{REMSisotropytheoremcorollaries} (1). We will now separate the cases where equality does or doesn't hold. \vspace{.5 \baselineskip}

\noindent {\it Case 1.} $\mydim{\phi_1} = 2^{n+1}$. In this case, \eqref{eq0} says that $\mydim{(\anispart{\phi_{K(\phi)_a})}} = 2^n$. Since $K(\phi)(\phi_1)_a = K(\phi)_a(\phi_1)$ is a purely transcendental extension of the function field of the quadric associated to $(\anispart{\phi_{K(\phi)_a})}$ (see \S \ref{SUBSECfunctionfields} above), it follows that $\mydim{(\anispart{\phi_{K(\phi)(\phi_1)_a})}} < 2^n$. Since the quasilinear quadratic form associated to $\anispart{(\mathfrak{b}_{K(\phi)(\phi_1)_a})}$ is a subform of $(\anispart{\phi_{K(\phi)(\phi_1)_a})}$, the Hauptsatz then gives that $\mathfrak{b}$ splits over $K(\phi)(\phi_1)_a$. Since $a$ is a value of $\pfister{X_1,\hdots,X_{n-2}}_{K(\phi)}$, however, $\pfister{X_1,\hdots,X_{n-2}}_{\mathfrak{b}}$, and hence $\mathfrak{c}$ is also split over $K(\phi)(\phi_1)_a$. Since $\mathfrak{b} = \anispart{(\mathfrak{c} \perp \mathfrak{d})}$, the same must then be true of $\mathfrak{d}$. By Theorem \ref{THMWittkernel}, we then also have that $a \in D((\phi_{\mathfrak{d}})_{K(\phi)(\phi_1)})$. Since $a \notin K(\phi)^2$, it follows that the isotropy index of $(2^n + 2^{n-2})$-dimensional quasiliner quadratic form
$$ \mathfrak{\eta}: = \pfister{X_1,\hdots,X_{n-2}} \perp  \phi_{\mathfrak{d}} $$
over $K(\phi)(\phi_1)$ is at least $2$ ($\pfister{a}$ is a common subform of the two summands). Observe, however that the isotropy index of $\eta$ over $K$ is just $1$ on account of $\lbrace X_1,\hdots,W_n \rbrace$ being a $2$-independent subset of $K$. Since $\mydim{\phi} > 2^n$, we then also have that $\witti{0}{\eta_{K(\phi)}} = 1$ by the separation theorem (Corollary \ref{CORcorofisotropytheorem} (4)), and so there exists a $(2^n + 2^{n-2} - 1)$-dimensional subform of $\phi$ which is anisotropic over $K(\phi)$ but becomes isotropic over $K(\phi)(\phi_1)$. By Corollary \ref{CORcorofisotropytheorem} (1), this implies that $\mydim{\phi_1} - \witti{1}{\phi_1} \leq 2^n + 2^{n-2} -2$. Since $\mydim{\phi_1} = 2^{n+1}$,  part (3) of the same result then implies that $\witti{1}{\phi_1} = 2^n = \frac{\mydim{\phi_1}}{2}$. By parts (5) and (6) of Corollary \ref{CORcorofisotropytheorem}, however, this implies that $\mathrm{lndeg}(\phi) = \mathrm{lndeg}(\phi_1) + 1 = n+2$, contradicting our earlier observation that $\phi$ is not a quasi-Pfister neighbour. The desired conclusion therefore holds in this case. \vspace{.5 \baselineskip}

\noindent {\it Case 2.} $\mydim{\phi_1} > 2^{n+1}$. In this case, the proof of Theorem \ref{THMgaptheoreminchar2} given above shows that there is an extension $L/K(\phi)$ with the following properties:

\begin{enumerate} \item[(i)] $L/K(\phi)$ is a tower of function fields of quadrics defined by anisotropic quasi-Pfister forms of foldness $\geq n+3$;
\item[(ii)] $(\phi_1)_L$ is an anisotropic quasi-Pfister neighbour. 
\end{enumerate}
By Proposition \ref{PROPcomplementarytoneighbours}, there then exists an anisotropic quasilinear quadratic form $\psi$ over $L$ such that $\mydim{\psi} = 2^{n+2} - \mydim{\phi_1}$ and $\witti{0}{(\phi_1)_{L(\psi)}} = \witti{1}{\psi} + \mydim{\phi_1} - 2^{n+1}$. Note that since $2^{n+1} < \mydim{\phi_1} < 2^{n+1} + 2^{n-1}$, we have $\mydim{\psi} > 2^n$. Since $\witti{0}{(\phi_1)_{L(\psi)}} = \witti{1}{\psi} + \mydim{\phi_1} - 2^{n+1}$, we have $\mydim{\anispart{(\phi_{L(\psi)})}} < 2^{n+1}$. As $\pfister{a}$ remains anisotropic over $L(\psi)$ (e.g., by the separation theorem), it is still a divisor of $\anispart{(\phi_{L(\psi)})}$ (Lemma \ref{LEMbehaviouroverfieldextensions} (2)), and so we have $\mydim{\anispart{(\phi_{L(\psi)_a})}}<2^n$. As already noted above, this implies (via the Hauptsatz) that $\mathfrak{b}_{L(\psi)_a}$ is split. Since we have already seen that $\mathfrak{c}$ splits over $K(\phi)_a$, and since $\mathfrak{b} = \anispart{(\mathfrak{c} \perp \mathfrak{d})}$, $\mathfrak{d}_{L(\psi)_a}$ must then also be split. Observe, however, that $\mathfrak{d}_{L_a}$ is not split. Indeed, consider the form $\eta$ of dimension $2^n + 2^{n-2}$ from Case 1 above. If $\mathfrak{d}$ were split over $L_a$, we would then have that $\witti{0}{\eta_L} \geq 2$ (for the same reasons as those given in Case 1). As before, this would imply that $\eta$ has a codimension-1 subform which is anisotropic over $K(\phi)$, but isotropic over $L$. In view of property (i) of $L$ however, the separation theorem (Corollary \ref{CORcorofisotropytheorem} (4)) rules out this possibility. Thus, $\mathfrak{d}_{L_a}$ is not split, and is hence anisotropic on account of being a Pfister form. Since $\mathfrak{d}$ splits over $L_a(\psi) = L(\psi)_a$, Theorem \ref{THMWittkernel} then tells us that $\mathrm{lndeg}(\psi_{L_a}) \leq n$, and so $\mathrm{lndeg}(\psi) = n+1$ by Corollary \ref{CORcorofisotropytheorem} (6). Since $\mydim{\psi} > 2^n$, it follows that $\psi$ is a quasi-Pfister neighbour, and so $\witti{1}{\psi} = \mydim{\psi} - 2^n = 2^{n+2} - \mydim{\phi_1} - 2^n$ by Remark \ref{REMSisotropytheoremcorollaries} (1). In particular, we have $\witti{0}{(\phi_1)_{L(\psi)}} = \witti{1}{\psi} + \mydim{\phi_1} - 2^{n+1} = 2^n$, and so $\mydim{\anispart{(\phi_{L(\psi)})}} = \mydim{\phi_1} -  2^n = 2^n + 2^{n-1} - 2x - \witti{1}{\phi}$. Now, the separation theorem again gives that $\mathfrak{c}$ remains anisotropic over $L$, and since $\mydim{\mathfrak{c}} < 2^{n+1}$, Theorem \ref{THMWittkernel} then tells us that $\mathfrak{c}$ does not split over $L(\psi)$. Since the quasilinear quadratic form associated to $(\mathfrak{c})_{L(\psi)}$ is a subform of $\anispart{((\phi_\mathfrak{c})_{L(\psi)})}$, the Hauptsatz then implies that the latter has dimension $\geq 2^n$. Note now that $\anispart{((\phi_\mathfrak{c})_{L(\psi)})}$ is divisible by the (still anisotropic!) quasi-Pfister form $\sigma: = \pfister{X_1,\hdots,X_{n-2}}$ by Lemma \ref{LEMbehaviouroverfieldextensions} (2). Since $\mydim{\anispart{((\phi_\mathfrak{c})_{L(\psi)})}} \geq 2^n$, it then follows from Corollary \ref{CORcorofisotropytheorem} (5) that the isotropy index of $\anispart{((\phi_\mathfrak{c})_{L(\psi)})}$ over $L(\psi)(\sigma)$ is at least $2^{n-1}$. As $\mathfrak{b}$ contains a codimension-$x$ subform of $\mathfrak{c}$, the isotropy index of $\anispart{(\phi_{L(\psi)})}$ over the same field is then at least $2^{n-1} -x$. But then $\mydim{\anispart{(\phi_{L(\psi)(\sigma)})}} \leq \mydim{\anispart{((\phi_\mathfrak{c})_{L(\psi)})}} - (2^{n-1} -x ) = 
2^n - x - \witti{1}{\phi}< 2^n$, which again forces $\mathfrak{b}_{L(\psi)(\sigma)}$ to split. Since $\mathfrak{c}$ also splits over $L(\psi)(\sigma)$ (because $\pfister{X_1,\hdots,X_{n-2}}_b$ does), and since $\mathfrak{b} = \anispart{(\mathfrak{c} \perp \mathfrak{d})}$, we then also get that the same is true of $\mathfrak{d}$. Now the separation theorem again tells us that $\mathfrak{d}$ remains anisotropic over $L$, and so Theorem \ref{THMWittkernel} then tells us that $\sigma_{L(\psi)}$ is a subform of $(\phi_{\mathfrak{d}})_{L(\psi)}$. Recalling that $\eta = \sigma \perp \phi_{\mathfrak{d}}$, we then get that $\anispart{(\eta_{L(\psi)})} \simeq (\phi_{\mathfrak{d}})_{L(\psi)}$, and so $\mathrm{lndeg}(\eta_{L(\psi)}) = n$. Note, however, that $\mathrm{lndeg}(\eta) = 2n - 2 \geq n+2$ by construction. In view of property (i) of $L$, Corollary \ref{CORcorofisotropytheorem} then implies that $\mathrm{lndeg}(\eta_{L(\psi)}) \geq n+1$, and so we have reached a contradiction. We are therefore also done in this case, and so the lemma is proved.  
\end{proof}
\end{lemma}

Putting Lemmas \ref{LEMadding2^n} and \ref{LEMtechnicaloptimalitylemma} together, we get:

\begin{proposition} \label{PROPmainoptimality} Let $m$ and $n$ be positive integers, and let $1 \leq i \leq n-1$. Let $0 \leq a, b \leq n-i+1$ be such that $a +b < n$, and set $r : = m2^{n-1} -2^{i-1} - 2^a - 2^b + 2$. Then $2^{n+1} + 2r \in \mathrm{Dim}(I^n, F)$. Moreover, if $n \geq 5$ and $i = n-1$, then $2^{n+1} + 2(r-1) \in \mathrm{Dim}(I^n,F)$. 
\begin{proof} By Lemma \ref{LEMadding2^n}, we can assume that $m=1$. By Lemma \ref{LEMtechnicaloptimalitylemma}, there then exist a field extension $K$ of $F$ and an anisotropic symmetric bilinear form $\mathfrak{b}$ of dimension $2^{n+1} + 2r$ over $K$ such that $[\mathfrak{b}] \in I^n(K)$. The first statement therefore holds. At the same time, in the case where $n\leq 5$ and $i = n-1$, Lemma \ref{LEMtechnicaloptimalitylemma} tells us that we can choose $\mathfrak{b}$ so that $\witti{1}{\phi_{\mathfrak{b}}} = 1$. By Lemma \ref{LEMWittindexrelation}, the Witt index of $\mathfrak{b}$ over $L: = F(\phi_{\mathfrak{b}})$ is then $1$, and so $\anispart{(\mathfrak{b}_L)}$ is an anisotropic symmetric bilinear form of dimension $2^{n+1} - 2(r-1)$ representing an element of $I^n(L)$, which proves the second statement.
\end{proof}
\end{proposition}

Now when $n \leq 5$, the reader will readily verify that every positive integer is expressible either in the form
\begin{itemize} \item $m2^{n-1} - 2^{i-1} - 2^a - 2^b +2$ for some positive integer $m$, integer $1 \leq i \leq n-1$ and integers $0 \leq a,b \leq n-i + 1$ with $a+b < n$, or
\item $m2^{n-1} - 2^{n-2} - 2^a - 2^b + 1$ for some positive integer $m$, and integers $0 \leq a,b \leq 2$. \end{itemize}
The previous proposition therefore gives the following:

\begin{corollary} \label{CORoptimalityfornleq5} If $n \leq 5$, then $\mathrm{Dim}(I^n,F)$ contains all even integers $\geq 2^{n+1}$. 
\end{corollary}

\begin{remarks} As far as Corollary \ref{CORoptimalityfornleq5} goes, the non-trivial part of Lemma \ref{LEMtechnicaloptimalitylemma} was only needed in the case where $n=5$, specifically to construct forms of dimension congruent to $2$ or $6$ modulo $32$. In the case where $n = 6$, one can expand the arguments from Lemma \ref{LEMtechnicaloptimalitylemma} in order to show that all even integers $\geq 128$ lie in $\mathrm{Dim}(I^6,F)$, with the possible exception of those which are now congruent to $2$ and $6$ modulo $64$. Since we were unable to resolve these cases, however, we refrain from giving the details here. In general, the main technical problem seems to be in producing examples of forms for which the norm degree of the associated quasilinear quadratic form is large enough in order to make splitting tower arguments feasible. We also note that Vishik's proof of the optimality of Theorem \ref{THMdimensiongaptheorem} in characteristic not 2 relies on the rigid behaviour exhibited by a certain discrete motivic invariant for low-dimensional forms representing elements of $I^n(F)$ that we alluded to in Remark \ref{REMSnonrigidity}. As we noted there, things are rather less rigid in the characteristic-2 setting, and so it does not seem impossible to the author that new dimension gaps beyond those captured by Theorem \ref{THMdimensiongaptheorem} may appear here. 
\end{remarks}

\section{Proof of Conjecture \ref{CONJAlbertconj} in Characteristic 2} \label{SECmain} Our goal now is to prove the characteristic-2 case of Conjecture \ref{CONJAlbertconj}. The key initial step is to show that if $\mathfrak{b}$ is an anisotropic symmetric bilinear form of dimension $2^n + 2^{n-1}$ over $F$ representing an element of $I^n(F)$, then its associated quasilinear quadratic form is divisible by an $(n-2)$-fold quasi-Pfister form over $F$. We will need the following simple specialization result, which may also be interpreted as the special case of \cite[Prop. 3.19]{Scully2} where $\mydim{\phi}$ is a multiple of $2^r$:

\begin{lemma} \label{LEMdivisibilityandpurelytranscendental} Let $\phi$ be an anisotropic quasilinear quadratic form over $F$, let $K$ be a purely transcendental extension of $F$, and let $r$ be a positive integer. If $\phi_K$ is divisible by an $r$-fold quasi-Pfister form, then the same is already true of $\phi$. 
\begin{proof} We can assume that $K = F(X)$ for a single interminate $X$. Suppose that $\pfister{f_1,\hdots,f_r}$ is an $r$-fold quasi-Pfister form over $K$ dividing $\phi$. Modifying by squares if needed, we can assume that $f_1,\hdots,f_r \in F[X]$. Let $a \in D(\phi)$. For each $1 \leq i \leq r$, we then have that $af_i \in D(\phi_K)$. By the Cassels-Pfister representation theorem (\cite[Thm. 17.3]{EKM}) and the additivity of $\phi$, however, any polynomial in $F[X]$ represented by $\phi$ lies in $D(\phi)[X^2]$ (see \cite[Cor. 3.4]{Hoffmann3}), so $f_i \in F[X^2]$ and $ab \in D(\phi)$ for every coefficient $b$ of $f_i$. The set $G(\phi)$ therefore contains all coefficients of the $f_i$. Since $\pfister{f_1,\hdots,f_r}$ is anisotropic over $K$, we must then have that $[G(\phi):F^2] \geq r$, and so $\mydim{\simform{\phi}}$ has foldness at least $r$, as desired. \end{proof} \end{lemma}

We will also need a couple of observations on the first higher isotropy indices of certain subforms of quasilinear quadratic forms. First, we have:

\begin{lemma} \label{LEMi1ofsubforms} Let $\phi$ and $\psi$ be anisotropic quasilinear quadratic forms of dimension $\geq 2$ over $F$ such that $\psi \subset \phi$. Suppose that the following hold:
\begin{enumerate} \item $\mathrm{lndeg}(\psi) = \mathrm{lndeg}(\phi) - 1$;
\item $\mathrm{lndeg}(\eta) = \mathrm{lndeg}(\phi)$ for every subform $\psi \subset \eta \subset \phi$ with $\eta \neq \psi$;
\item $\mydim{\psi} \geq \frac{\mydim{\phi} + \witti{1}{\psi}}{2}$. \end{enumerate}
Then $\witti{1}{\psi} \geq \witti{1}{\phi}$. 
\begin{proof} We may assume that $1 \in D(\psi)$. Let $\sigma \subset \phi$ be such that $\phi \simeq \psi \perp \sigma$. Since $\mydim{\psi} \geq \frac{\mydim{\phi} + \witti{1}{\psi}}{2}$, we have $\mydim{\sigma} = r$ for some $r \leq \mydim{\psi} - \witti{1}{\psi}$. Now, as $\mathrm{lndeg}(\psi) = \mathrm{lndeg}(\phi) - 1$, there exists an element $x \in F^\times$ such that $N(\phi) = N(\psi)(x)$. We can therefore find elements $a_1,\hdots,a_r,b_1,\hdots,b_r \in N(\psi)$ such that $\sigma \simeq \form{a_1 + xb_1,\hdots, a_r + xb_r}$. Observe that the form $\form{b_1,\hdots,b_r}$ is anisotropic. Indeed, if this were not the case, there would exist scalars $\lambda_1,\hdots,\lambda_r \in F$, at least one of which is non-zero, such that $\sum_{i=1}^r \lambda_i^2 b_i = 0$. Since $\sigma$ is anisotropic, it would then follow that
$$ y: = \sum_{i=1}^r \lambda_i^2a_i = \left( \sum_{i=1}^r \lambda_i^2a_i\right)+ x\left(\sum_{i=1}^r \lambda_i^2 b_i\right) = \sum_{i=1}^r \lambda_i^2(a_i + xb_i) $$
is a non-zero element of $D(\sigma)$. But $y \in N(\psi)$, and so $\eta: = \psi \perp \form{y}$ would then be a subform of $\phi$ with $\mathrm{lndeg}(\eta) = \mathrm{lndeg}(\phi) - 1$, contradicting (2). Since $r \leq \mathrm{dim}(\psi) - \witti{1}{\psi}$, Corollary \ref{CORcorofisotropytheorem} (1) then tells us that $\form{b_1,\hdots,b_r}$ remains anisotropic over $F(\psi)$. Now, let $V$ be the underlying vector space of $\phi$. We may assume that $\psi$ is the restriction of $\phi$ to a subspace $U \subset V$. We claim that all vectors in $V\otimes_F F(\psi)$ which are isotropic with respect to $\phi$ lie in $U \otimes_F F(\psi)$. Suppose, for the sake of contradiction, that this is not the case. There then exists a vector $u \in U \otimes_F F(\psi)$ and elements $\lambda_1,\hdots,\lambda_r \in F(\psi)$, at least one of which is non-zero, such that $\psi(u) + \sum_{i=1}^r \lambda_i^2(a_i + xb_i) = 0$. Set $z: = \sum_{i=1}^r \lambda_i^2 b_i \in N(\psi)$. Since at least one of the $\lambda_i$ is non-zero, and since $\form{b_1,\hdots,b_r}$ remains anisotropic over $F(\psi)$, $z$ is non-zero. Rearranging the preceeding identity, we then get that $x = z^{-1}(\psi(u) + \sum_{i=1}^r \lambda_i^2 a_i)$. Since $1 \in D(\psi)$, however, we have $\psi(u) \in N(\psi)$, and so this shows that $x \in N(\psi_{F(\psi)}) = N(\psi_1)$. Then $N(\phi_{F(\psi)}) =  N(\psi_1)(x) = N(\psi_1)$, and so $\mathrm{lndeg}(\phi_{F(\psi)}) = \mathrm{lndeg}(\psi_1) = \mathrm{lndeg}(\psi) - 1 = \mathrm{lndeg}(\phi) - 2$. As this contradicts Corollary \ref{CORcorofisotropytheorem} (6), our initial claim must then hold, and so $\witti{0}{\phi_{F(\psi)}} = \witti{1}{\psi}$. Since $\witti{1}{\phi}$ is the smallest non-zero isotropy index attained by $\phi$ over all extensions of $F$, this proves the lemma.
\end{proof}
\end{lemma}

A variant of this that we will also use is the following:

\begin{lemma} \label{LEMi1ofsubforms2} Let $\phi$ be an anisotropic quasilinear quadratic form of dimension $\geq 2$ such that $\phi \subset x\pi \perp \tau$ for some anisotropic quasilinear quadratic form $\tau$ over $F$ with $1 \in D(\tau)$, scalar $x \in F \setminus N(\tau)$ and quasi-Pfister form $\pi \subset \normform{\tau}$. Let $\psi$ be the subform of $\phi$ with value set $D(\phi) \cap D(\tau)$, and assume that $\mydim{\psi} \geq 2$. If $\witti{1}{\psi} < \witti{1}{\phi}$, then $\normform{\psi} \subset \pi$. 
\begin{proof} Let $\sigma \subset \phi$ be such that $\phi \simeq \psi \perp \sigma$, and set $r: = \mydim{\sigma}$. Since $\phi \subset \tau \perp x\pi$, we have $\sigma \simeq \form{a_1 + xb_1,\hdots,a_r + xb_r}$ for some $a_1,\hdots,a_r \in D(\tau)$ and $b_1,\hdots,b_r \in D(\pi)$. By our choice of $\psi$, we have $D(\sigma) \cap D(\tau) = 0$. Arguing as in the previous lemma, it then follows that $\form{b_1,\hdots,b_r}$ is anisotropic, and hence a subform of $\pi$. Now, since $\pi \subseteq \normform{\tau}$, and since $1 \in D(\tau)$, the elements $a_1,\hdots,a_r, b_1,\hdots,b_r$ all lie in $N(\tau)$. Since $x$ does not lie in $N(\tau)$, the arguments from the proof of the previous lemma then show that $\witti{1}{\psi} \geq \witti{1}{\phi}$ unless $\form{b_1,\hdots,b_r}$ becomes isotropic over $F(\psi)$. But $\form{b_1,\hdots,b_r} \subset \pi$, so if the latter occurs, then $\normform{\psi} \subset \pi$ by Corollary \ref{CORcorofisotropytheorem} (6).
\end{proof}
\end{lemma}

Now, let $\phi$ be an anisotropic quasilinear quadratic form of dimension $2^n$ for some integer $n \geq 2$. By Remark \ref{REMSisotropytheoremcorollaries} (1), the largest possible value of the integer $\witti{1}{\phi}$ is $2^{n-1}$, and this is realized precisely when $\phi$ is similar to a quasi-Pfister form by Corollary \ref{CORcorofisotropytheorem} (5). By part (3) of the same corollary, the next largest possible value is $2^{n-2}$. Here, Theorem \ref{THMkeytool} allows us to prove the following:

\begin{theorem} \label{THMi12^{n-2}} Let $n$ be an integer $\geq 2$, and let $\phi$ be an anisotropic quasilinear quadratic form of dimension $2^n$ over $F$. Then $\witti{1}{\phi} = 2^{n-2}$ if and only if $\mathrm{lndeg}(\phi) = n+1$ and $\phi$ is divisible by an $(n-2)$-fold quasi-Pfister form over $F$.
\begin{proof} Suppose first that the stated conditions hold. Since $\mathrm{lndeg}(\phi) > n$, $\phi$ is not similar to a quasi-Pfister form, and so $\witti{1}{\phi} \leq 2^{n-2}$ by the remarks preceding the statement of the theorem. To prove that equality holds, it suffices to show that $\mydim{\phi_1} = \mydim{\phi} - \witti{1}{\phi}$ is divisible by $2^{n-2}$. But $\phi$ is divisible by an $(n-2)$-fold quasi-Pfister form $\pi$ over $F$, and this form remains anisotropic over $F(\phi)$ by the separation theorem (Corollary \ref{CORcorofisotropytheorem} (4)). By Lemma \ref{LEMbehaviouroverfieldextensions} (2), $\phi_1$ is then divisible by $\pi_{F(\phi)}$, and so the claim holds. 

Conversely, suppose that $\witti{1}{\phi} = 2^{n-2}$. By Corollary \ref{CORcorofisotropytheorem} (2), $\phi_1$ is then an anisotropic form of dimension $2^{n-1} + 2^{n-2}$ which is divisible by a quasi-Pfister form of foldness at least $n-2$. We must therefore have that $\phi_1 \simeq \nu \otimes \form{x,y,z}$ for some $(n-2)$-fold quasi-Pfister form $\nu$ over $F(\phi)$ and scalars $x,y,z \in F(\phi)^\times$. By Remark \ref{REMSisotropytheoremcorollaries} (2), we then have that $\mathrm{lndeg}(\phi) = \mathrm{lndeg}(\phi_1) + 1 \leq n +1$, and since $\phi$ is not a quasi-Pfister form (because $\witti{1}{\phi} \neq 2^{n-1}$), equality must in fact hold. It now remains to show that $\phi$ is divisible by an $(n-2)$-fold quasi-Pfister form. We claim here the following:

\begin{claim} \label{CLAIMPfisterneighboursubform} In the above situation, there exists a purely transcendental extension $K$ of $F$ such that $\phi_K$ admits a subform of dimension $> 2^{n-1}$ which is a quasi-Pfister neighbour. 
\end{claim}

Before proving the claim, let us first show how it gives what we want. By Lemma \ref{LEMdivisibilityandpurelytranscendental}, we are free to pass to purely transcendental extensions of $F$. Given the claim, we can therefore assume that $\phi$ admits a subform of dimension $> 2^{n-1}$ which is a quasi-Pfister neighbour, i.e., has norm degree $n$. Among all such subforms, choose one of largest possible dimension, say $\psi$. Since $\mydim{\psi} > 2^{n-1} = \frac{\mydim{\phi}}{2}$, Lemma \ref{LEMi1ofsubforms} tells us that $\witti{1}{\psi} \geq \witti{1}{\phi} = 2^{n-2}$. By Corollary \ref{CORcorofisotropytheorem} (3), this implies that $\mydim{\psi} \geq 2^{n-1} + 2^{n-2}$. On the other hand, Corollary \ref{CORcorofisotropytheorem} (6) tells us that $\psi$ remains anisotropic over $F(\phi)$. By Lemma \ref{LEMbehaviouroverfieldextensions} (1), we then have that $\psi_{F(\phi)} \subset \phi_1$, and so $\mydim{\psi} \leq 2^{n-1} + 2^{n-2}$. Thus, $\psi$ has dimension exactly $2^{n-1} + 2^{n-2}$.  By Proposition \ref{PROPcomplementarytoneighbours}, there then exists a quasilinear quadratic form $\sigma$ of dimension $2^{n-2}$ over $F$ such that $\witti{0}{\psi}_{F(\sigma)} > 2^{n-2}$. Now we have already noted that $\psi_{F(\phi)}$ is a subform of $\phi_1$, and since both forms have dimension $2^{n-1} + 2^{n-2}$, they are in fact isometric. In particular, we have $\witti{0}{\phi_{F(\phi)(\sigma)}} = 2^{n-2} + \witti{0}{\psi}_{F(\phi)(\sigma)} \geq 2^{n-2} + \witti{0}{\psi_{F(\sigma)}} > 2^{n-1}$. We claim that $\witti{0}{\phi_{F(\sigma)}} \geq 2^{n-1}$. In other words, we claim that $\eta: = \anispart{(\phi_{F(\sigma)})}$ has dimension $\leq 2^{n-1}$. Suppose not. By Corollary \ref{CORcorofisotropytheorem} (3), we then have that $\mydim{\eta_1} \geq 2^{n-1}$, so that $\witti{0}{\phi_{F(\sigma)(\eta)}} \leq 2^{n-1}$. But $F(\sigma)(\eta)$ is a purely transcendental extension of $F(\sigma)(\phi) = F(\phi)(\sigma)$ (see \S \ref{SUBSECfunctionfields} above), so this contradicts our earlier observation that $\witti{0}{\phi_{F(\phi)(\sigma)}} > 2^{n-1}$. We thus have that $\witti{0}{\phi_{F(\sigma)}} \geq 2^{n-1}$, and parts (1) and (5) of Corollary \ref{CORcorofisotropytheorem} then tell us that $\phi$ is divisible by $\normform{\sigma}$. Since $\mydim{\sigma} = 2^{n-2}$, $\normform{\sigma}$ has foldness at least $n-2$, and so $\phi$ is divisible by an $(n-2)$-fold quasi-Pfister form, as desired. We now prove Claim \ref{CLAIMPfisterneighboursubform}. 

Recall here that $\phi_1 \simeq \nu \otimes \form{x,y,z}$ for some $(n-2)$-fold quasi-Pfister form $\nu$ over $F(\phi)$ and scalars $x,y,z \in F(\phi)^\times$. Per \S \ref{SUBSECfunctionfields}, the field $F(\phi)$ is $F$-isomorphic to an inseparable quadratic extension of a purely transcendental extension of $F$, say $K_a$ with $K/F$ purely transcendental and $a \in K \setminus K^2$. By Lemma \ref{LEMquadraticextensions}, (1), we have $D(\phi_1) \subseteq K$, so $x,y,z \in K^\times$, and we may also view $\nu$ as a quasi-Pfister form over $K$. Consider the form $\pi := \pfister{a} \otimes \nu$ over $K$. By Lemma \ref{LEMquadraticextensions} (2), $\pi \otimes \form{x,y,z}$ is an anisotropic form of dimension $2^n + 2^{n-1}$ over $K$ whose value set coincides with $D(\phi_1)$. Since $D(\phi_K) \subset D(\phi_1)$, $\phi_K$ is then a subform of $\pi \otimes \form{x,y,z}$. Set $\tau: = \pi \otimes \form{y,z}$, so that $\phi_K \subset x\pi \perp \tau$. Replacing $\phi_K$ with $y\phi_K$ if needed (this doesn't affect our goal), we can assume that $y = 1$, and hence that $\tau$ is an $n$-fold quasi-Pfister form containing $\pi$ as a subform. Since $\mathrm{lndeg}(\phi) = n+1$, we then also have that $x \notin N(\tau)$. Let $\psi$ now be the subform of $\phi_K$ whose value set is $D(\phi_K) \cap D(\tau)$. Since $\mydim{\pi} = 2^{n-1}$, $\psi$ has dimension at least $2^{n-1}$. If $\mydim{\psi} > 2^{n-1}$, then $\psi$ is a quasi-Pfister neighbour and we're done. Suppose now that $\mydim{\psi} = 2^{n-1}$. If $\normform{\psi} \subset \pi$, then $\mathrm{lndeg}(\psi) \leq \mathrm{lndeg}(\pi) = n-1$, and so $\psi$ is similar to a quasi-Pfister form. If $\normform{\psi} \not \subset \pi$, then Lemma \ref{LEMi1ofsubforms2} tells us that $\witti{1}{\psi} \geq \witti{1}{\phi_K} = 2^{n-2}$. By the remarks preceding the statement of the theorem, this again means that $\psi$ is similar to a quasi-Pfister form. Thus, $\psi$ is always similar to a quasi-Pfister form, and any $(2^n + 1)$-dimensional subform of $\phi_K$ that contains it as a subform is therefore a quasi-Pfister neighbour . This completes the proof of Claim \ref{CLAIMPfisterneighboursubform}, and hence the proof of the theorem.
\end{proof}
\end{theorem}

\begin{remark} \label{REMi12^{n-2}} If $q$ is a non-singular quadratic form of dimension $2^n \geq 4$ over a field (of any characteristic) with $\witti{1}{q} = 2^{n-2}$, then it is expected that $q$ is divisible by an $(n-2)$-fold Pfister form. At present, however, this is only known to be true for $n \leq 4$. Theorem \ref{THMi12^{n-2}} therefore also solves the quasilinear analogue of a non-trivial open problem in the theory of non-singular quadratic forms. In each case, the known restrictions on the invariant $\mathfrak{i}_1$ (see Remark \ref{REMSisotropytheoremcorollaries} (1)) allow us to interpret the problem as one of classifying the anisotropic forms for which the value of $\mathfrak{i}_1$ is exactly a quarter of the dimension of the form. \end{remark}

Now, to achieve the initial goal, we shall also need the following theorem, which is the quasilinear version of a result of Vishik in the theory of non-singular quadratic forms:

\begin{theorem}[{\cite[Thm. 9.2]{Scully1}}] \label{THMVishikcomparison} Let $\phi$ be an anisotropic quasilinear quadratic form over $F$ with $\mydim{\phi} = 2^n + m$ for some positive integer $n$ and integer $1 \leq m \leq 2^n$. If $L$ is a field extension of $F$, then either $\witti{0}{\phi_L} \geq m$ or $\witti{0}{\phi_L} \leq m - \witti{1}{\phi}$. 
\end{theorem} 

We can now prove:

\begin{proposition} \label{PROPdivisibilityinmainresult} Let $\mathfrak{b}$ be an anisotropic symmetric bilinear form over $F$ of dimension $2^n + 2^{n-1}$ for some integer $n \geq 2$, and let $\phi = \phi_{\mathfrak{b}}$. If $[b] \in I^n(F)$, then one of the following holds:
\begin{enumerate} \item $\mathrm{lndeg}(\phi) = n+1$ and $\phi$ is divisible by an $(n-1)$-fold quasi-Pfister form over $F$;
\item $\mathrm{lndeg}(\phi) = n+2$ and $\phi$ is divisible by an $(n-2)$-fold quasi-Pfister form over $F$.
\end{enumerate}
\begin{proof} Set $\mathfrak{b}_1: = \anispart{(\mathfrak{b}_{F(\phi)})}$. Since $\mathfrak{b}$ becomes isotropic over $F(\phi)$, $\mydim{\mathfrak{b}_1} < 2^{n} + 2^{n-1}$. Since $[\mathfrak{b}_1] \in I^n(F(\phi))$, Theorem \ref{THMgaptheoreminchar2} then tells us that $\mydim{\mathfrak{b}_1} \leq 2^n$. If this equality were strict, the same result would tell us that $\mathfrak{b}_1$ is trivial, i.e., that $\mathfrak{b}$ splits over $F(\phi)$. Since $\mydim{\normform{\phi}} \geq 2^{n+1}$, however, this possibility is eliminated by Theorem \ref{THMWittkernel}. We therefore have that $\mydim{\mathfrak{b}_1} = 2^n$, and so $\mathfrak{b}_1$ is a general $n$-fold bilinear Pfister form by Theorem \ref{THMHauptsatz}. In particular, $\phi_{\mathfrak{b}_1} = y\eta$ for some $y \in F(\phi)^\times$ and anisotropic quasi-Pfister form $\eta$ over $F$. Now the Witt index of $\mathfrak{b}_{F(\phi)}$ is $\frac{\mydim{\mathfrak{b}} - \mydim{\mathfrak{b}_1}}{2} = 2^{n-2}$, and so $\witti{1}{\phi}$ is at least $2^{n-2}$ by Lemma \ref{LEMWittindexrelation}. If $\witti{1}{\phi} > 2^{n-2}$, then Corollary \ref{CORcorofisotropytheorem} (3) tells us that $\witti{1}{\phi} = 2^{n-1}$, and so $\mydim{\phi_1} = 2^n$. But $\phi_{\mathfrak{b}_1}$ is a subform of $\phi_1$, and so $\phi_1 \simeq \phi_{\mathfrak{b}_1}$ in this case. In particular, $\phi_1$ is similar to a quasi-Pfister form, and so $\mathrm{lndeg}(\phi) = \mathrm{lndeg}(\phi_1) + 1 = n+1$ by Remark \ref{REMSisotropytheoremcorollaries} (2). In other words, $\phi$ is a quasi-Pfister neighbour. The proof of Theorem \ref{THMgaptheoreminchar2} then shows that $\phi$ is divisible by an $(n-1)$-fold quasi-Pfister form. More precisely, we can choose here a bilinear lift of $\phi$ which is embedded as a subform of an $(n+1)$-fold bilinear Pfister form. If we let $\mathfrak{c}$ be the complementary subform of this bilinear Pfister form, then the proof of Theorem \ref{THMgaptheoreminchar2} shows that $\phi_{\mathfrak{c}}$ is an $(n-1)$-quasi-Pfister form that divides $\phi$. Suppose now that $\witti{1}{\phi} = 2^{n-2}$. By Corollary \ref{CORcorofisotropytheorem} (2), $\phi_1$ is then a form of dimension $2^n + 2^{n-2}$ which is divisible by an $(n-2)$-fold quasi-Pfister form, say $\nu$. By Corollary \ref{CORcorofisotropytheorem} (5), we then have that $\witti{0}{(\phi_1)_{F(\phi)(\nu)}} = 2^{n-1} + 2^{n-3}$. Since $\phi_{\mathfrak{b}_1} \simeq y\eta$ is a codimension-$2^{n-2}$ subform of $\phi_1$, it must then become isotropic over $F(\phi)(\sigma)$. By Corollary \ref{CORcorofisotropytheorem} (6), we then have that $\nu \subset \eta$. Let $x$ be an element of $D(\phi_1)$ which does not lie in $D(\phi_{\mathfrak{b}_1}) = yD(\eta)$. Since $\phi_1$ is divisible by $\nu$, $x\nu$ is a subform of $\phi_1$. Since $D(\eta)$ is a subfield of $F$ containing $D(\nu)$, we have $xD(\nu) \cap yD(\eta) = 0$. As a result, $x \nu \perp y\eta$ is anisotropic and hence a subform of $\phi_1$. Since the two forms have the same dimension, we then have that $\phi_1 \simeq x\nu \perp y\eta$. Since $\nu \subset \eta$ (but $x \notin yD(\eta)$), this shows that $\mathrm{lndeg}(\phi_1) = n+1$, and so $\mathrm{lndeg}(\phi) = n+2$. It remains to show that $\phi$ is divisible by an $(n-2)$-fold quasi-Pfister form. As in the proof of Theorem \ref{THMi12^{n-2}}, we use that $F(\phi)$ is $F$-isomorphic to an inseparable quadratic extension of a purely transcendental extension of $F$, say $K_a$ with $K/F$ purely transcendental and $a \in K \setminus K^2$. As before, the elements $x,y$ then lie in $K$, and we can view $\eta$ and $\nu$ as quasi-Pfister forms over $K$. Set $\tau: = y\pfister{a} \otimes \eta$ and $\pi : = \pfister{a} \otimes \nu$. Then $x\pi \perp \tau$ is an anisotropic form of dimension $2^{n+1} + 2^{n-1}$ over $K$ containing $\phi_K$ as a subform (again, see the proof of Theorem \ref{THMi12^{n-2}}). We are going to show that $\phi_K$ is divisible by an $(n-2)$-fold quasi-Pfister form. To this end, we can replace $\phi_K$ with $y\phi_K$ and therefore assume that $y = 1$. Then $\tau$ is an $(n+1)$-fold quasi-Pfister form containing $\pi$ as as a subform. Moreover, since $\mathrm{lndeg}(\phi_K) = \mathrm{lndeg}(\phi) = n+2 > n+1 = \mathrm{lndeg}(\tau)$, we have $x \notin N(\tau)$. Let $\psi$ be the subform of $\phi_K$ with value set $D(\phi_K) \cap D(\tau)$. Since $\mydim{\pi} = 2^{n-1}$, we have $\mydim{\psi} \geq \mydim{\phi} - 2^{n-1} = 2^n$. In particular, $\normform{\psi}$ has dimension at least $2^n$, and is therefore not a subform of $\pi$. By Lemma \ref{LEMi1ofsubforms2}, it follows that $\witti{1}{\psi} \geq \witti{1}{\phi_K} = 2^{n-2}$. We claim that $\mydim{\psi} = 2^n$. Suppose, to the contrary, that this is not the case. Since $\tau$ is an $(n+1)$-fold quasi-Pfister form, $\psi$ is then a quasi-Pfister neighbour. Since $\phi_K$ is not a quasi-Pfister neighbour, we then have that $\mydim{\psi} < \mydim{\phi} = 2^n + 2^{n-1}$. On the other hand, since $\witti{1}{\psi} \geq 2^{n-2}$, Corollary \ref{CORcorofisotropytheorem} (3) tells us that $\mydim{\psi} \geq 2^n + 2^{n-2}$. We therefore have that $2^n + 2^{n-2} \leq \mydim{\psi} < 2^n + 2^{n-1}$. Now, since $\psi$ is a quasi-Pfister neighbour, Proposition \ref{PROPcomplementarytoneighbours} tells us that there exists an anisotropic quasilinear quadratic form $\rho$ over $K$ such that $\mydim{\rho} = 2^{n+1} - \mydim{\psi} > 2^{n-1}$ and $\witti{0}{\psi_{K(\rho)}} > \mydim{\psi} - 2^n \geq 2^{n-2}$. In particular, we have $\witti{0}{\phi_{K(\rho)}} > 2^{n-2}$. Since $\mydim{\phi} = 2^n + 2^{n-1}$ and $\witti{1}{\phi_K} = 2^{n-2}$, Theorem \ref{THMVishikcomparison} then tells us that $\witti{0}{\phi_{K(\rho)}} \geq 2^{n-1}$. In other words, $\anispart{(\phi_{K(\rho)})}$ has dimension at most $2^n$. Now, since $\mydim{\rho} > 2^{n-1}$, we have $\mathrm{lndeg}(\rho) \geq n$. Since $\mydim{\mathfrak{b}}$ is not divisible by $2^n$, it then follows from Theorem \ref{THMWittkernel} that $\mathfrak{b}$ does not split over $K(\rho)$. Since $\mathfrak{b}$ becomes isotropic over this field, Theorem \ref{THMHauptsatz} then tells us that $\anispart{(\mathfrak{b}_{K(\rho)})}$ is similar to a Pfister form. Since the quasilinear quadratic form associated to $\anispart{(\mathfrak{b}_{K(\rho)})}$  is a subform of $\anispart{(\phi_{K(\rho)})}$, and since the latter has dimension at most $2^n$, we conclude that $\anispart{(\phi_{K(\rho)})}$ is similar to an $n$-fold quasi-Pfister form. In particular, $\mathrm{lndeg}(\anispart{(\phi_{K(\rho)})}) = n$. But Corollary \ref{CORcorofisotropytheorem} (6), then implies that $\mathrm{lndeg}(\phi_K) = n+1$, a contradiction. We must therefore have that $\mydim{\psi} = 2^n$. Now $\witti{1}{\psi} \geq 2^{n-2}$, so Corollary \ref{CORcorofisotropytheorem} (3) tells us that either $\witti{1}{\psi} = 2^{n-2}$ or $\witti{1}{\psi} = 2^{n-1}$. In the first case, Theorem \ref{THMi12^{n-2}} tells us that $\psi$ is divisible by an $(n-2)$-fold quasi-Pfister form. In the second, $\psi$ is similar to a quasi-Pfister form (see the remarks preceding Theorem \ref{THMi12^{n-2}}). Thus, in either case, $\psi$ is divisible by an $(n-2)$-fold quasi-Pfister form, say $\alpha$. We claim that $\phi_K$ is divisible by $\alpha$. To see this, note that we have $\witti{0}{\phi_{K(\alpha)}} \geq \witti{0}{\psi_{K(\alpha)}} = \frac{\mydim{\psi}}{2} = 2^{n-1}$ by Corollary \ref{CORcorofisotropytheorem} (5), and so $\anispart{(\phi_{K(\alpha)})}$ has dimension at most $2^n$. If the bilinear form $\mathfrak{b}$ did not split over $K(\alpha)$, it would then follow that $\mathrm{lndeg}(\anispart{(\phi_{K(\alpha)})}) = n$, and hence that $\mathrm{lndeg}(\phi_K) = n+1$ (this is identical to the argument used to show that $\mydim{\psi} = 2^n$ above). Since this is not the case, $\mathfrak{b}_{K(\alpha)}$ must be split, and Theorem \ref{THMWittkernel} then implies that $\phi_K$ is then divisible by $\alpha$, as desired. But since $K$ is a purely transcendental extension of $F$, Lemma \ref{LEMdivisibilityandpurelytranscendental} then tells us that $\phi$ is already divisible by an $(n-2)$-fold quasi-Pfister form over $F$, thereby completing the proof. \end{proof} \end{proposition}

To deduce Theorem \ref{THMmain} from this, we shall need two additional results. The first is well known, and essentially due to Elman and Lam, who considered the analogous statement in characteristic not 2 (see \cite[Prop. 24.5]{EKM}). For the sake of completeness, we quickly explain how the characteristic-2 case follows from Theorem \ref{THMWittkernel}:

\begin{lemma} \label{LEMlinkageof3Pfisters} Let $n$ be an integer $\geq 2$, and let $\mathfrak{b}_1$, $\mathfrak{b}_2$ and $\mathfrak{b}_3$ be anisotropic $(n-1)$-fold bilinear Pfister forms over $F$. If $[\mathfrak{b}_1] + [\mathfrak{b}_2] + [\mathfrak{b}_3] \in I^n(F)$, then there exists an $(n-2)$-fold bilinear Pfister form $\mathfrak{c}$ over $F$ and elements $a_1,a_2\in F^\times$ such that $\mathfrak{b}_1 \simeq  \mathfrak{c} \otimes \pfister{a_1}_b$, $\mathfrak{b}_2 \simeq \mathfrak{c} \otimes \pfister{a_2}_b$ and $\mathfrak{b}_3 \simeq \mathfrak{c} \otimes \pfister{a_1a_2}_b$. 
\begin{proof} For each $i$, set $\phi_i: = \phi_{\mathfrak{b}_i}$. Consider now the form $\mathfrak{b} := \anispart{(\mathfrak{b}_1 \perp \mathfrak{b}_2)}$, which represents an element of $I^{n-1}(F)$. Since $1 \in D(\phi_1) \cap D(\phi_2)$, $\mathfrak{b}_1 \perp \mathfrak{b}_2$ is isotropic, and so $\mydim{\mathfrak{b}} < 2^n$. If $\mathfrak{b}$ were trivial, we would have that $[\mathfrak{b}_3] \in I^n(F)$, contradicting the Hauptsatz. Since $[\mathfrak{b}] \in I^{n-1}(F)$, the latter therefore tells us that $\mydim{\mathfrak{b}} \geq 2^{n-1}$. On the other hand, $\mathfrak{b}_3$ splits over $F(\phi_3)$, and so $[\mathfrak{b}_{F(\phi_3)}] \in I^n(F(\phi_3))$. Since $\mydim{\mathfrak{b}} < 2^n$, another application of the Hauptsatz then gives that $\mathfrak{b}$ splits over $F(\phi_3)$. By Theorem \ref{THMWittkernel}, it follows that $\mathfrak{b} \simeq \lambda \mathfrak{d}$ for some scalar $\lambda \in F^\times$ and anisotropic $(n-1)$-fold quasi-Pfister form $\mathfrak{d}$ over $F$ with $\phi_{\mathfrak{d}} \simeq \phi_3$. Since $\lambda \mathfrak{d}$ and $\mathfrak{d}$ have the same class in $W(F)/I^n(F)$, we then have that $[\mathfrak{d} \perp \mathfrak{b}_3] = [\mathfrak{d}] + [\mathfrak{b}_3] = [\mathfrak{b}_1] + [\mathfrak{b}_2] + [\mathfrak{b}_3]\in I^n(F)$. But since $\phi_{\mathfrak{d}} \simeq \phi_3$, Lemma \ref{LEMliftsofquasiPfisters} tells us that the anisotropic part of $\mathfrak{d} \perp \mathfrak{b}_3$ has dimension at most $2^{n-1}$. By the Hauptsatz, it follows that $\mathfrak{d} \perp \mathfrak{b}_3$ is split, and so $\mathfrak{b}_3 \simeq \mathfrak{d}$. Finally, since $\mydim{\mathfrak{b}} = 2^{n-1}$, the linkage theorem for bilinear Pfister forms (see \cite[Proposition 6.21]{EKM}) tells us that there exists an $(n-2)$-fold Pfister form $\mathfrak{c}$ over $F$ and scalars $a_1,a_2 \in F^\times$ such that $\mathfrak{b}_1 \simeq  \mathfrak{c} \otimes \pfister{a_1}_b$ and $\mathfrak{b}_2 \simeq \mathfrak{c} \otimes \pfister{a_2}_b$. Then $\mathfrak{b} \simeq \mathfrak{c} \otimes \form{a_1,a_2}_{b} $, and so $\mathfrak{b}_3 \simeq \mathfrak{d} \simeq \mathfrak{c} \otimes \pfister{a_1a_2}_b$. 
\end{proof}
\end{lemma}

Recall now that if $\pi$ is an anisotropic quasi-Pfister form over $F$, then we write $\mathscr{A}(\pi)$ for the set consisting of the bilinear Pfister forms over $F$ whose associated quadratic form is isometric to $\pi$. The second result we need is:

\begin{proposition} \label{PROPdeterminantlemma} Let $n$ be an integer $\geq 2$, let $\pi = \pfister{a_1,\hdots,a_{n-2}}$ be an anisotropic $(n-2)$-fold quasi-Pfister form over $F$, and let $\mathfrak{b}_1,\hdots,\mathfrak{b}_r \in \mathscr{A}(\pi)$. Suppose that $x_2,\hdots,x_{r-1} \in F^\times$ are such that $\lbrace a_1,\hdots,a_{n-2},x_2,\hdots,x_{r-1} \rbrace$ is a $2$-independent subset of $F$, and let $x_r \in F^\times$. If $[\mathfrak{b}_1 \perp x_2\mathfrak{b}_2 \perp \cdots \perp x_r\mathfrak{b}_r] \in I^n(F)$, then $r$ is even, $\mathfrak{b}_1 \simeq \cdots \simeq \mathfrak{b}_r$, and $x_r \in x_2\cdots x_{r-1}F^2(a_1,\hdots,a_{n-2})$. 
\begin{proof} We use Kato's Theorem \ref{THMKato}. Let $2 \leq i \leq r$. Since $\mathfrak{b}_i \in \mathscr{A}(\pi)$, we have $\mathfrak{b}_i = \pfister{\alpha_1,\hdots,\alpha_{n-2}}_b$ for some $\alpha_1,\hdots,\alpha_{n-2} \in D(\pi) = F^2(a_1,\hdots,a_{n-2})$. In $\Omega_F^{n-2}$, the wedge product $\frac{d\alpha_1}{\alpha_1} \wedge \cdots \wedge \frac{d\alpha_{n-2}}{\alpha_{n-2}}$ is then equal to $\mu_i \frac{da_1}{a_1} \wedge \cdots \wedge \frac{da_{n-2}}{a_{n-2}}$ for some $\mu_i \in F^2(a_1,\hdots,a_{n-2})$. On the other hand, this wedge product is by definition the image of the class of $\mathfrak{b}_i$ under the injective homomorphism $e_{n-2} \colon I^{n-2}/I^{n-1}(F) \rightarrow \Omega_F^{n-2}$. Since $\mathfrak{b}_i$ is anisotropic, its class in $I^{n-2}(F)/I^{n-1}(F)$ is nontrivial by the Hauptsatz, and so we must have that $\mu_i \neq 0$. Consider now the form $\mathfrak{b}: = \mathfrak{b}_2 \otimes \pfister{x_2}_b \perp \cdots \perp \mathfrak{b}_r \otimes  \pfister{x_{r-1}}_b $. Since $\mathfrak{b}_2,\hdots,\mathfrak{b}_r$ are $(n-2)$-fold Pfister forms, we have $[\mathfrak{b}] \in I^{n-1}(F)$. Since $[\mathfrak{b}] = [\mathfrak{b}_1 \perp x_2\mathfrak{b}_2 \perp \cdots \perp x_r\mathfrak{b}_r] + [\mathfrak{b}_1 \perp \mathfrak{b}_2 \perp \cdots \perp \mathfrak{b}_r]$ in $W(F)$, and since $[\mathfrak{b}_1 \perp x_2\mathfrak{b}_2 \perp \cdots \perp x_r\mathfrak{b}_r] \in I^n(F)$, it follows that $[\mathfrak{b}_1 \perp \mathfrak{b}_2 \perp \cdots \perp \mathfrak{b}_r] \in I^{n-1}(F)$. By Lemma \ref{LEMliftsofquasiPfisters}, however, the anisotropic part of $\mathfrak{b}_1 \perp \mathfrak{b}_2 \perp  \cdots \perp \mathfrak{b}_r$ has dimension at most $2^{n-2}$. By the Hauptsatz, it follows that $[\mathfrak{b}_1 \perp \mathfrak{b}_2 \perp \cdots \perp \mathfrak{b}_r] = 0$, and so $\mathfrak{b} = [\mathfrak{b}_1 \perp x_2\mathfrak{b}_2 \perp \cdots \perp x_r\mathfrak{b}_r] \in I^n(F)$. In particular, the homomorphism $e_{n-1} \colon I^{n-1}(F)/I^n(F) \rightarrow \Omega_F^{n-1}$ annihilates the class of $\mathfrak{b}$. Since $e_{n-1}([\mathfrak{b}_i \otimes \pfister{x_i}_b]) = e_{n-2}([\mathfrak{b}_i]) \wedge \frac{dx_i}{x_i}$ for all $i$, it follows that
\begin{equation} \label{eq1} \sum_{i=2}^r \mu_i \frac{da_1}{a_1} \wedge \cdots \wedge \frac{da_{n-2}}{a_{n-2}}\wedge \frac{dx_i}{x_i} = 0 \end{equation}
in $\Omega_F^{n-1}$. Since $\lbrace a_1,\hdots,a_{n-2},x_2,\hdots,x_{r-1} \rbrace$ is a $2$-independent subset of $F$, the first $r-2$ terms in the sum on the left side of this equation are $F$-linearly independent. Moreover, in order for the equality to hold, the set $\lbrace a_1,\hdots,a_{n-2},x_2,\hdots,x_{r-1},x_r \rbrace$ must be $2$-dependent, i.e., we have $x_r \in F^2(a_1,\hdots,a_{n-2},x_2,\hdots,x_{r-1})$. For each $2 \leq i \leq r-1$, set $F_i: = F^2(a_2,\hdots,a_{n-2},x_2,\hdots,\widehat{x_i},\hdots,x_{r-1})$, where the hat indicates omission of the element $x_i$. By the $2$-independence of the set $\lbrace a_1,\hdots,a_{n-2},x_2,\hdots,x_{r-1} \rbrace$, $F^2(a_1,\hdots,a_{n-2},x_2,\hdots,x_{r-1})$ is the internal direct sum of $F_i$ and $x_iF_i$. In particular, we have $x_r = u + x_iv$ for unique elements $u,v \in F_i$. We claim that $u = 0$. To see this, note that we have $dx_r = du + x_idv + vdx_i$ in $\Omega_F$. Since $u,v \in F_i$, both $du$ and $dv$ can be written as $F$-linear combinations of the elements $da_1,\hdots,da_{n-2},dx_2,\hdots,\widehat{dx_i},\hdots,dx_{r-1} \in \Omega_F$. If we substitute our expression for $dx_r$ into the left side of equation \eqref{eq1}, we will therefore obtain an $F$-linear combination of the elements $\frac{da_1}{a_1} \wedge \cdots \wedge \frac{da_{n-2}}{a_{n-2}}\wedge \frac{dx_j}{x_j}$ with $2 \leq j \leq r-1$. Since these elements are $F$-linearly independent, the coefficient of $\frac{da_1}{a_1} \wedge \cdots \wedge \frac{da_{n-2}}{a_{n-2}}\wedge \frac{dx_i}{x_i}$ must be $0$. From the preceding discussion, however, this coefficient is nothing else but $\mu_i + x_i\mu_rvx_r^{-1} $. Since $\mu_i,\mu_r \in F^2(a_1,\hdots,a_{n-2}) \subseteq F_i$, and since $\mu_i \neq 0$, we then have that $x_r = x_i(\mu_i^{-1}\mu_{r}v) \in x_iF_i $, and so our claim holds. Note that we also get that $\mu_{i}^{-1}\mu_rv = v$, and so $\mu_i = \mu_r$ on account of $x_r$ being non-zero. In particular, the class of the form $\mathfrak{b}_i \perp \mathfrak{b}_r$ is annihilated by the homomorphism $e_{n-2}$, and so $[\mathfrak{b}_i \perp \mathfrak{b}_r] \in I^{n-1}(F)$ by the injectivity of the latter. Again, however, a combination of Lemma \ref{LEMliftsofquasiPfisters} and the Hauptsatz then tells us that $\mathfrak{b}_i \perp \mathfrak{b}_r$ is split, and hence that $\mathfrak{b}_i \simeq \mathfrak{b}_r$ by the anisotropy of $\mathfrak{b}_i$ and $\mathfrak{b}_r$. This shows that $\mathfrak{b}_2 \simeq \cdots \simeq \mathfrak{b}_r$. At the same time, we have already seen that $\mathfrak{b}_1 \perp \mathfrak{b}_2 \perp \cdots \perp \mathfrak{b}_r$ is split, so this also gives that $[\mathfrak{b}_1] = (r-1)\cdot [\mathfrak{b}_r]$ in $W(F)$. Since $\mathfrak{b}_1$ is not split, it then follows that $r$ is even, and that $[\mathfrak{b}_1] = [\mathfrak{b}_r]$. Again, by the anisotropy of $\mathfrak{b}_1$ and $\mathfrak{b}_r$, this gives that $\mathfrak{b}_1 \simeq \mathfrak{b}_r$, so $r$ is even and $\mathfrak{b}_1 \simeq \cdots \simeq \mathfrak{b}_r$. It now only remains to show that $x_r \in x_2 \cdots x_{r-1}F^2(a_1,\hdots,a_{n-2})$. But since $\lbrace a_1,\hdots,a_{n-2},x_2,\hdots,x_{r-1} \rbrace$ is a $2$-independent subset of $F$, the set $x_2 \cdots x_{r-1}F^2(a_1,\hdots,a_{n-2})$ is precisely the intersection of the sets $x_iF_i$, and so the needed assertion also follows from the preceding discussion. This completes the proof of the proposition. \end{proof} \end{proposition}

We are now ready to prove Theorem \ref{THMmain}. For the reader's convenience, we restate it here, along with the converse, which holds trivially.

\begin{theorem} Let $n$ be an integer $\geq 2$, and $\mathfrak{b}$ be an anisotropic symmetric bilinear form of dimension $2^n + 2^{n-1}$ over $F$. Then $[\mathfrak{b}] \in I^n(F)$ if and only if $\mathfrak{b} \simeq \mathfrak{c} \otimes \mathfrak{d}$ for some $(n-2)$-fold bilinear Pfister form $\mathfrak{c}$ and Albert form $\mathfrak{d}$ over $F$. 
\begin{proof} If $\mathfrak{d}$ is an Albert form over $F$, then $[\mathfrak{d}] \in I^2(F)$, and so the if part is clear. Conversely, suppose that $[\mathfrak{b}] \in I^n(F)$. Set $\phi: = \phi_{\mathfrak{b}}$. We now separate two cases. \vspace{.5 \baselineskip}

\noindent {\it Case 1.} $\phi$ is a quasi-Pfister neighbour, i.e., $\mathrm{lndeg}(\phi) = n+1$. In this case, Proposition \ref{PROPdivisibilityinmainresult} tells us that $\phi$ is then divisible by an $(n-1)$-fold quasi-Pfister form over $F$, say $\pi$. By Corollary \ref{CORcorofisotropytheorem} (5), we then have that $\witti{0}{\phi_{F(\pi)}} = \frac{\mydim{\phi}}{2}$, and so $\mydim{\anispart{(\phi_{F(\pi)})}} = \frac{2^n + 2^{n-1}}{2} < 2^n$. Since the quasilinear quadratic form associated to $\anispart{(\mathfrak{b}_{F(\pi)})}$ is a subform of $\anispart{(\phi_{F(\pi)})}$, we then have that $\mydim{\anispart{(\mathfrak{b}_{F(\pi)})}} < 2^n$, and so $\mathfrak{b}$ splits over $F(\pi)$ by the Hauptsatz. By Theorem \ref{THMWittkernel}, we then have that $\mathfrak{b} \simeq x_1\mathfrak{b}_1 \perp x_2 \mathfrak{b}_2 \perp x_3 \mathfrak{b}_3$ for some $x_1,x_2,x_3 \in F^\times$ and $(n-1)$-fold bilinear Pfister forms $\mathfrak{b}_1,\mathfrak{b}_2,\mathfrak{b}_3 \in \mathscr{A}(\eta)$. Now multiplying an $(n-1)$-fold bilinear Pfister form over $F$ by a non-zero scalar does not change its class in $W(F)/I^n(F)$, and so $[\mathfrak{b}_1] + [\mathfrak{b}_2] + [\mathfrak{b}_3] = [x_1\mathfrak{b}_1] + [x_2\mathfrak{b}_2] + [x_3\mathfrak{b}_3] = [\mathfrak{b}] \in I^n(F)$. By Lemma \ref{LEMlinkageof3Pfisters}, there then exist an $(n-2)$-fold bilinear Pfister form $\mathfrak{c}$ over $F$ and scalars $a_1,a_2 \in F^\times$ such
that $\mathfrak{b}_1 \simeq \mathfrak{c} \otimes \pfister{a_1}_b$, $\mathfrak{b}_2 \simeq \mathfrak{c} \otimes \pfister{a_2}_b$ and $\mathfrak{b}_3 \simeq \mathfrak{c} \otimes \pfister{a_1a_2}_b$. Factoring, we get that $\mathfrak{b} \simeq \mathfrak{c} \otimes \mathfrak{d}$, where $\mathfrak{d} = \form{x_1,x_1a_1,x_2,x_2a_2,x_3,x_3a_1a_2}_b$. Since the latter has trivial determinant, this proves the theorem in this case. \vspace{.5 \baselineskip}

\noindent {\it Case 2.} $\phi$ is not a quasi-Pfister neighbour. In this case, Proposition \ref{PROPdivisibilityinmainresult} tells us that $\mathrm{lndeg}(\phi) = n+2$, and that $\phi$ is divisible by an $(n-2)$-fold quasi-Pfister form, say $\pi = \pfister{a_1,\hdots,a_{n-2}}$. As in Case 1, Corollary \ref{CORcorofisotropytheorem} (5) then tells us that $\mydim{\anispart{(\phi_{F(\pi)})}}< 2^n$, which forces $\mathfrak{b}$ to split over $F(\pi)$ by the Hauptsatz. By Theorem \ref{THMWittkernel}, there then exist scalars $x_1,\hdots,x_6 \in F^\times$ and $(n-2)$-fold bilinear Pfister forms $\mathfrak{b}_1,\hdots,\mathfrak{b}_6 \in \mathscr{A}(\pi)$ such that $\mathfrak{b} \simeq x_1\mathfrak{b}_1 \perp \cdots  \perp x_6\mathfrak{b}_6$. Replacing $\mathfrak{b}$ with $x_1^{-1}\mathfrak{b}$, we can assume that $x_1 = 1$. Then $N(\phi) = F^2(a_1,\hdots,a_{n-2},x_2,\hdots,x_6)$. Since $\mathrm{lndeg}(\phi) = n+2$, we can assume (after possibly reordering the $x_i$) that $\lbrace a_1,\hdots,a_{n-2},x_2,\hdots,x_5 \rbrace$ is a $2$-independent subset of $F$. Since $[\mathfrak{b}] \in I^n(F)$, Proposition \ref{PROPdeterminantlemma} then tells us that $\mathfrak{b}_1 \simeq \cdots \simeq \mathfrak{b}_6$, and that $x_6 \in x_2\hdots x_5 F^2(a_1,\hdots,a_{n-2})$, say $x_6 = x_2\hdots x_5 \lambda$ with $\lambda \in F^2(a_1,\hdots,a_{n-2})$. But since $\mathfrak{b}_6 \in \mathscr{A}(\pi)$, and since $D(\pi) = F^2(a_1,\hdots,a_{n-2})$, the roundness of Pfister forms gives that $\lambda \mathfrak{b}_6 \simeq \mathfrak{b}_6$. Replacing $x_6$ with $x_6\lambda^{-1}$, we can therefore assume that $x_6 = x_2\cdots x_5$. But then $\mathfrak{b} \simeq \mathfrak{c} \otimes \mathfrak{d}$ with $\mathfrak{c} := \mathfrak{b}_1$ being an $(n-2)$-fold bilinear Pfister form and $\mathfrak{d} := \form{1,x_2,\hdots,x_6}_b$ being an Albert form. This completes this case, and so the theorem is proved.  \end{proof} \end{theorem}
\vspace{.5 \baselineskip}

\noindent {\bf Acknowledgements.} This work was supported by NSERC Discovery Grant No. RGPIN-2019-05607. The author is grateful to Detlev Hoffmann for comments on an initial draft. 

\bibliographystyle{alphaurl}

\end{document}